\documentclass[10pt,twoside]{article}
\usepackage{mathrsfs}
\usepackage{amsmath}
\usepackage{amssymb}
\usepackage{fancyhdr}
\usepackage{latexsym}
\usepackage{bbding}
\usepackage{mathrsfs}
\usepackage{wasysym}
\usepackage{cite}

\usepackage{multicol,graphics}

\newtheorem{theorem}{Theorem}[section]
\newtheorem{lemma}[theorem]{Lemma}
\newtheorem{corollary}[theorem]{Corollary}

\newtheorem{definition}[theorem]{Definition}
\newtheorem{remark}[theorem]{Remark}

\numberwithin{equation}{section}
\newenvironment{proof}[1][Proof]{\noindent\textbf{#1.} }{\hfill $\Box$}
\allowdisplaybreaks

 \makeatletter
\begin{document}
\title{{Gevrey analyticity of solutions to the 3D nematic liquid crystal flows in critical Besov space}
\thanks{
This work is partially supported by the National Natural Science Foundation of China (11401202), the Scientific Research Fund of Hunan Provincial
Education Department (14B117), and the China Postdoctoral Science Foundation (2015M570053).
}
}
\author{
   {\small  
   Qiao Liu$^{1,2}$ \thanks{\text{Corresponding author. E-mail address}: liuqao2005@163.com. }}
\\
{\small  $^{1}$\textit{Department of Mathematics, Hunan Normal University, Changsha, Hunan 410081, China}}\\
{\small $^{2}$\textit{Institute of Applied Physics and Computational Mathematics, Beijing, 100088,  China}}
}
\date{}
\maketitle

\begin{abstract}

We show that the solution to the Cauchy problem of the 3D nematic
 liquid crystal flows, with initial data
 belongs to a critical Besov space, belongs to a Gevrey class.
More precisely, it is proved that for any 
 $(u_{0},d_{0}-\overline{d}_{0})\in \dot{B}^{\frac{3}{p}-1}_{p,1}(\mathbb{R}^{3})\times
\dot{B}^{\frac{3}{q}}_{q,1}(\mathbb{R}^{3})$ with some suitable conditions imposed on $p, q\in(1,\infty)$, 
there exists $T^{*}>0$ depending only on initial
data, such that the nematic liquid crystal flows admits a unique solution $(u,d)$
on $\mathbb{R}^{3}\times (0,T^{*})$, and satisfies
\begin{align*}
\|e^{\sqrt{t}\Lambda_{1}}{u}(t)\|_{\widetilde{L}^{\infty}_{T^{*}}(\dot{B}^{\frac{3}{p}-1}_{p,1})\cap
\widetilde{L}^{1}_{T^{*}}(\dot{B}^{\frac{3}{p}+1}_{p,1})}
+\|e^{\sqrt{t}\Lambda_{1}}({d}(t)-\overline{d}_{0})\|_{\widetilde{L}^{\infty}_{T^{*}}(\dot{B}^{\frac{3}{q}}_{q,1})
\cap \widetilde{L}^{1}_{T^{*}}(\dot{B}^{\frac{3}{q}+2}_{q,1})}< \infty.
\end{align*}
Here, $\overline{{d}}_{0}\in \mathbb{S}^{2}$ is a constant unit
vector, and $\Lambda_{1}$ is the Fourier multiplier whose symbol is given by $|\xi|_{{1}}=|\xi_{1}|+|\xi_{2}|+|\xi_{3}|$.
Moreover, if the initial data is sufficiently small enough, then $T^{*}=\infty$.
As a consequence of the method, decay estimates of higher-order derivatives of solutions
in Besov spaces are deduced.


\end{abstract}

\section{Introduction}\label{Int}

\noindent

Liquid crystal, which is a state of matter capable of flow, but its
molecules may be oriented in a crystal-like way. Research into liquid crystals is an area of a very successful synergy between mathematics and
physics (see  \cite{IS}).
To our knowledge, three main phases of liquid crystals are
distinguished, nematic, termed smectic and cholesteric. The nematic
phase appears to be the most common one, where the molecules do not
exhibit any positional order, but they have long-range orientational
order. 
The Ericksen-Leslie system is one of the most successful models for the nematic
liquid crystals. It was formulated by Ericksen and  Leslie in
1960s (see \cite{ER,LE}), who derived suitable constitutive equations.
 Since then, many remarkable developments have been made from both theoretical and applied aspects.

In the present paper, we investigate  the Cauchy problem of the following three dimensional
simplified version of the Ericksen-Leslie system  in the whole
space, which describes the motion of the incompressible flow of
nematic liquid crystals:
\begin{align}\label{eq1.1}
\begin{cases}
{\partial_{t}}u-\Delta u +(u\cdot\nabla)u+\nabla{P}=-\nabla\cdot(\nabla d \odot\nabla d),\quad&(x,t)\in\mathbb{R}^{3}\times (0,+\infty),\\
 \partial_{t}d+(u\cdot\nabla)d=\Delta d+|\nabla
d|^{2}d,\quad\quad
  \quad\quad\quad\quad\quad\quad&(x,t)\in\mathbb{R}^{3}\times (0,+\infty),\\
\quad\nabla\cdot u=0,\quad\quad\quad\quad\quad\quad\quad\quad
\quad\quad\quad\qquad\quad\quad\; &(x,t)\in\mathbb{R}^{3}\times (0,+\infty),\\
(u,d)|_{t=0}=(u_{0},d_{0}),\quad\quad\quad\quad\quad\quad\quad\quad\quad\quad\quad\quad
&x\in\mathbb{R}^{3},
\end{cases}
\end{align}
where $u(x,t):\mathbb{R}^{3}\times (0,+\infty)\rightarrow
\mathbb{R}^{3}$ is the unknown velocity field of the flow,
$P(x,t):\mathbb{R}^{3}\times (0,+\infty)\rightarrow
\mathbb{R}$ is the scalar pressure,
$d(x,t):\mathbb{R}^{3}\times (0,+\infty)\rightarrow \mathbb{S}^{2}$,
the unit sphere in $\mathbb{R}^{3}$, is the unknown (averaged)
macroscopic/continuum molecule orientation of the nematic liquid
crystal flow,  $\nabla \cdot u=0$ represents
the incompressible condition, $u_{0}$ is a given initial velocity
with $\nabla \cdot u_{0}=0$ in distribution sense,  and $d_{0}:
\mathbb{R}^{3}\rightarrow \mathbb{S}^{2}$ is a given initial liquid
crystal orientation field and satisfies
 $\lim_{|x|\rightarrow \infty} d_{0}(x)=\overline{d}_{0}$
 with the constant unit vector $\overline{d}_{0}\in \mathbb{S}^{2}$. The notation
$\nabla d\odot\nabla d$ denotes the $3\times 3$ matrix whose
$(i,j)$-th entry is given by $\partial_{i}d\cdot
\partial_{j}d$ ($1\leq i,j\leq 3$), and there holds
$\nabla\cdot(\nabla d\odot \nabla d)= \Delta d\cdot\nabla
d+\frac{1}{2}\nabla |\nabla d|^{2}$.
Throughout the paper, due to the concrete values of  the viscosity constants do not play a special
role in our discussion, we assume that they  are all equal
to one for simplicity.

The above simplified
 Ericksen--Leslie model \eqref{eq1.1} was first introduced by Lin \cite{L}.
  When $d\equiv \overline{d}_{0}$,
it corresponds to the following well-known Navier--Stokes equations
\begin{align}
    \label{eq1.2}
\begin{cases}
{\partial_{t}}u-\Delta u +(u\cdot\nabla)u+\nabla{P}=-0,\quad &(x,t)\in\mathbb{R}^{3}\times (0,+\infty),\\
\quad \nabla\cdot u=0,\quad &(x,t)\in\mathbb{R}^{3}\times (0,+\infty),\\
u|_{t=0}=u_{0}, &x\in\mathbb{R}^{3}.
\end{cases}
\end{align}
This equations has drawn much attention of researchers, and has
been extensively studied in the past several decades.
 Notice that \eqref{eq1.2} is scaling invariant in the following sense: if $(u,P)$ solves \eqref{eq1.2},
 so does $(u_{\lambda},P_{\lambda}):=(\lambda u(\lambda x,\lambda^{2}t), \lambda^{2}P(\lambda x,\lambda^{2}t))$
  with initial
data $u_{0\lambda}:=\lambda u_{0}(\lambda x)$. We say a function space $X$ defined in
$\mathbb{R}^{3}$, which is said to be the initial critical space for \eqref{eq1.2}, if
the associated norm is invariant under the transformation $u_{0}\rightarrow u_{0\lambda}$
    (up to a positive constant independent of $\lambda$).                    
Under these scalings, it is easy to see that
$L^{3}(\mathbb{R}^{3})$, $\dot{H}^{\frac{1}{2}}(\mathbb{R}^{3})$, $\dot{B}^{\frac{3}{p}-1}_{p,q}(\mathbb{R}^{3})$
and $BMO^{-1}$ are initial critical spaces for \eqref{eq1.2}, and one can find various well-posedness results for the initial data
belongs to these initial critical spaces in \cite{MC,TK,KTataru,PG} and the references therein.
 For the issue of regularity of solutions,  Foias-Temam \cite{FT} developed an energy method, which involves certain  pseudo-differential operators,
to prove that the solutions of \eqref{eq1.2}, starting with initial data in the Sobolev space $H^{1}(\mathbb{R}^{3})$,
 become instantaneously elements
of a certain Gevrey class of regularity, which, in
particular, makes them real analytic functions. Using iterative derivative
estimates, the analyticity of solutions to \eqref{eq1.2} for  small initial data in
$L^{3}(\mathbb{R}^{3})$ was obtained by Giga-Sawada \cite{GS}, and in $BMO^{-1}$ was obtained by
Germain-Pavlovic-Staffilani \cite{GPS} (see also Miura and Sawada
\cite{MS} on the iterative derivative techniques).
 Very recently, base on the method introduced by Foias-Temam \cite{FT},
 Bae-Biswas-Tadmor \cite{BBT} and Huang-Wang \cite{CHBW}
established analyticity of solutions to  \eqref{eq1.2} through the Gevery estimate in Besov space.
 More precisely, they proved that for initial data $u_{0}\in \dot{B}^{-1+\frac{3}{p}}_{p,q}(\mathbb{R}^{3})$
 with $1<p<\infty$, $1\leq q\leq \infty$, there exists an unique solution $u(t)$ defined
  on $ \mathbb{R}^{3}\times(0,T)$ such that
\begin{align*}
\|e^{\sqrt{t}\Lambda_{1}} u\|_{\widetilde{L}^{\infty}_{T}(\dot{B}^{-1+\frac{3}{p}}_{p,q})
\cap \widetilde{L}^{1}_{T}(\dot{B}^{1+\frac{3}{p}}_{p,q}) }\leq C\|u_{0}\|_{\dot{B}^{-1+\frac{3}{p}}_{p,q}},
\end{align*}
where $\Lambda_{1}$ is the Fourier multiplier whose symbol is given by $|\xi|_{1}=|\xi_{1}|+|\xi_{2}|+|\xi_{3}|$.

For system \eqref{eq1.1}, due to the appearance of the nonlinear term $|\nabla d|^{2}d$ with
 the restriction $|d|=1$ causes significant mathematical difficulties, Lin-Liu
\cite{LL1,LL2} have initiated the mathematical analysis of \eqref{eq1.1} by considering its Ginzburg--Landau
approximation. Namely, the Dirichlet energy $\int_{\mathbb{R}^{3}}\frac{1}{2}|\nabla d|^{2}\text{d}x$ for
 $d:\mathbb{R}^{3}\rightarrow \mathbb{S}^{2}$ is replaced by
the Ginzburg-Landau energy $\int_{\mathbb{R}^{3}}(\frac{1}{2}|\nabla d|^{2}
+\frac{1}{4\varepsilon^{2}}(1-|d|^{2})^{2})\text{d}x$ $(\varepsilon>0)$
for $d: \mathbb{R}^{3}\rightarrow \mathbb{R}^{3}$,  that is, equation $\eqref{eq1.1}_{2}$
is replaced by
\begin{align*}
d_{t}+u\cdot\nabla d=\Delta d+ \frac{1}{\varepsilon^{2}}(1-|d|^{2})d.
\end{align*}
 In  this situation, Lin-Liu in
\cite{LL1} proved the local existence of classical solutions and the
global existence of weak solutions in dimensions two and three with
Dirichlet boundary conditions, and in \cite{LL2} established the existence of suitable weak solutions
and their partial regularity in dimensions three, analogous to the celebrated regularity results by
Caffarelli-Kohn-Nirenberg \cite{CKN} for the three dimensional Navier--Stokes equations.
For more studies about the Ginzburg--Landau
approximation system, we refer the readers to
\cite{HW,XLW} and the references therein.

In the past a few years, progress has also been made on the analysis of system \eqref{eq1.1}.
Lin-Lin-Wang \cite{LLW}
established that there exists global Leray--Hopf type weak solutions
to the initial boundary value problem for system
\eqref{eq1.1} on bounded domains in two space
dimensions (see also \cite{HMC}). The uniqueness of such weak
solutions is proved by Lin-Wang \cite{LW}, see also Xu-Zhang \cite{XZ} for related works.
Very recently, when the space dimension is three, Lin-Wang \cite{LW14}
established the existence of global weak solutions when the initial
 data $(u_{0},d_{0})\in L^{2}\times H^{1}$ with the initial director field $d_{0}$ maps to
 the upper hemisphere $\mathbb{S}_{+}^{2}$,
 yet the global existence of weak solutions to system \eqref{eq1.1}
  with general initial data in dimensions three is still not resolved.
  As for an issue of the existence of strong solutions, Wen and Ding
\cite{WD} obtained local existence and uniqueness of strong
solution, Hineman and Wang \cite{JHW} established the global
well--posedness of system \eqref{eq1.1} in dimensions
three with small initial data $(u_{0},d_{0})$ in $L^{3}_{uloc}$,
where $L^{3}_{uloc}$ is the space of uniformly locally
$L^{3}$-integrable functions in $\mathbb{R}^{3}$,
Wang \cite{W} proved the global-in-time existence of strong
solutions for the incompressible liquid crystal model in the whole space provided that the initial data $({u}_{0},{d}_{0})$
\begin{align*}
\|{u}_{0}\|_{BMO^{-1}}+[{d}_{0}]_{BMO}<\varepsilon,
\end{align*}
for some suitable small positive $\varepsilon$. For the issue of regularity of solutions,
Du and Wang \cite{DW}, and Liu
\cite{QL1} proved that the small global strong solution obtained by
\cite{W} is arbitrary space-time regularity and is algebraically
decay as time goes to infinity. On the other hand, in order to understand
which quantity goes to infinite as the time approaches to infinity,
various blow-up criteria have been established in \cite{LLW,HW1} and the
references therein.

Similar as the Navier--Stokes equations, system
\eqref{eq1.1} is invariant under the following
transformation
\begin{align*}
(u_{\lambda}(x,t),P_{\lambda}(x,t),d_{\lambda}(x,t)):= (\lambda u(\lambda x, \lambda^{2}t),
\lambda^{2}P(\lambda x,\lambda^{2}t), d(\lambda x, \lambda^{2}t)).
\end{align*}
We say a functional space is the initial critical
   space for system \eqref{eq1.1} if the associated norm is invariant under the transformation
   $(u_{0},d_{0})\rightarrow (u_{0\lambda}, d_{0\lambda}):= (\lambda u_{0}(\lambda x), d_{0}(\lambda x))$
    is invariant for all $\lambda >0$. In fact, the homogeneous Besov space
    $\dot{B}^{\frac{3}{p}-1}_{p,1}(\mathbb{R}^{3})\times \dot{B}^{\frac{3}{q}}_{q,1}(\mathbb{R}^{3})$
    with $1<p,q<\infty$
     is the initial critical space for system \eqref{eq1.1}.
In this paper, motivated by the papers \cite{AB,BBT,FT,CHBW,PG},
we focus on regularity of solutions to the nematic liquid crystal flows \eqref{eq1.1}. In particular, we establish that
the solutions to  system
\eqref{eq1.1}, with initial data
$(u_{0},d_{0}-\overline{d}_{0})$ belongs to the initial critical Besov spaces
$\dot{B}^{\frac{3}{p}-1}_{p,1}(\mathbb{R}^{3})\times \dot{B}^{\frac{3}{q}}_{q,1}(\mathbb{R}^{3})$
with some suitable conditions imposed on  $p, q\in(1,\infty)$, immediately becomes Gevrey regular on
the time existence interval. In order to state our main result,
let us first define the Gevrey regular on $L^{p}$-based Besov spaces.

\begin{definition}\label{def1.1}
(See \cite{BBT,AB}) We say that a function $f$ is \textit{Gevrey regular} if
\begin{align*}
\|e^{\gamma\Lambda} f\|_{\dot{B}^{s}_{p,q}}<\infty,
\end{align*}
for some $s\in\mathbb{R}, \gamma>0$, and $1\leq p,q\leq \infty$.
 Here, the \textit{Gevrey operator } $e^{\gamma\Lambda}$ is the  Fourier  multiplier operator whose
 symbol is given by $e^{\gamma |\xi|}$,
 where $|\xi|=\sqrt{|\xi_{1}|^{2}+|\xi_{2}|^{2}+|\xi_{3}|^{2}}$.
\end{definition}

Our main result is as follows:

\begin{theorem}\label{thm1.1}
Let $1<p,q<\infty$ such that
\begin{align}\label{eq1.3}
-\inf\{\frac{1}{3},\frac{1}{2p}\}\leq \frac{1}{q}-\frac{1}{p}\leq \frac{1}{3}.
\end{align}
Let $\overline{d}_{0}\in \mathbb{S}^{2}$ be a constant vector,
$u_{0}\in\dot{B}^{\frac{3}{p}-1}_{p,1}(\mathbb{R}^{3})$ with
$\operatorname{div} u_{0}=0$,
$d_{0}-\overline{d}_{0}\in\dot{B}^{\frac{3}{q}}_{q,1}(\mathbb{R}^{3})$, and
\begin{align*}
M_{0}:=\|u_{0}\|_{\dot{B}^{\frac{3}{p}-1}_{p,1}}+\|d_{0}-\overline{d}_{0}\|_{\dot{B}^{\frac{3}{q}}_{q,1}}.
\end{align*}
Then there exists a $T^{*}= T^{*}(M_{0})>0$ such that system
\eqref{eq1.1} has a unique solution  $(u,d)$ on $\mathbb{R}^{3}\times (0,T^{*})$, and
\begin{align*}
\!(u,d\!-\!\overline{d}_{0})\!\in
&\!\widetilde{L}^{\!\infty}\!(0,\!T^{*}; e^{\!\sqrt{t}\Lambda_{1}}\!\dot{B}^{\frac{3}{p}\!-\!1}_{p,1}\!
(\mathbb{R}^{\!3})\!)\!
\cap \!\widetilde{L}^{1}\!(0,\!T^{*}; e^{\!\sqrt{t}\Lambda_{1}}\!
\dot{B}^{\frac{3}{p}\!+\!1}_{p,1}\!(\mathbb{R}^{\!3})\!)
\!\times\! \widetilde{L}^{\!\infty}\!(0,\!T^{*}; e^{\!\sqrt{t}\Lambda_{1}}\!\dot{B}^{\frac{3}{q}}_{q,1}\!
(\mathbb{R}^{\!3})\!)\!
\cap \!\widetilde{L}^{1}\!(0,\!T^{*}; e^{\!\sqrt{t}\Lambda_{1}}\dot{B}^{\frac{3}{q}\!+\!2}_{q,1}\!(\mathbb{R}^{\!3})\!).
\end{align*}
If $T^{*}<\infty$, then we have for some $\theta\in (0,1]$
\begin{align}\label{eq1.4}
&\|u\|_{\widetilde{L}^{1}_{T^{*}}( e^{\sqrt{t}\Lambda_{1}}\dot{B}^{\frac{3}{p}+1}_{p,1})}
+\|u\|_{\widetilde{L}^{1+\theta}_{T^{*}}
(e^{\sqrt{t}\Lambda_{1}}\dot{B}^{\frac{3}{p}+\frac{1-\theta}{1+\theta}}_{p,1})\cap
\widetilde{L}^{\frac{1+\theta}{\theta}}_{T^{*}}
(e^{\sqrt{t}\Lambda_{1}}\dot{B}^{\frac{3}{p}+\frac{\theta-1}{1+\theta}}_{p,1})}\nonumber\\
&+\|d-\overline{d}_{0}\|_{\widetilde{L}^{1}_{T^{*}}(  e^{\sqrt{t}\Lambda_{1}}\dot{B}^{\frac{3}{q}+2}_{q,1})}
+ \|d-\overline{d}_{0}\|_{\widetilde{L}^{1+\theta}_{T^{*}}
(e^{\sqrt{t}\Lambda_{1}}\dot{B}^{\frac{3}{q}+\frac{2}{1+\theta}}_{q,1})\cap
\widetilde{L}^{\frac{1+\theta}{\theta}}_{T^{*}}
(e^{\sqrt{t}\Lambda_{1}}\dot{B}^{\frac{3}{q}+\frac{2\theta}{1+\theta}}_{q,1})}
=\infty.
\end{align}
Moreover, if $M_{0}$ is sufficiently small, then we have $T^{*}=\infty$.

\end{theorem}

\begin{remark}\label{rem1.3}
1.\ In Definition \ref{def1.1}, $\dot{B}^{s}_{p,q}$ denotes the homogeneous $L^{p}$-based Besov spaces with
regularity index $s$ and summability index $q$ (see Section 2 below), we notice that when $p=q=2$, one recovers
the usual definition of Gevrey classes (cf. \cite{FT,PG}) in the space-periodic setting.
\medskip

2.\ We say $f\in
\widetilde{L}^{\rho}(0,T; e^{\sqrt{t}\Lambda_{1}}\dot{B}^{s}_{p,r}(\mathbb{R}^{\!3}))$,  
$s\in\mathbb{R}$, $1\leq \rho,r\leq \infty$,
if and only if $ e^{\sqrt{t}\Lambda_{1}} f\in
\widetilde{L}^{\rho}(0,T;\dot{B}^{s}_{p,r}(\mathbb{R}^{\!3}))$. Moreover,  it holds that
$\|f\|_{\widetilde{L}^{\rho}_{T}(e^{\sqrt{t}\Lambda_{1}}\dot{B}^{s}_{p,r})}
=\|e^{\sqrt{t}\Lambda_{1}} f\|_{\widetilde{L}^{\rho}_{T}(\dot{B}^{s}_{p,r})}$.
\medskip

3.\ We emphasize that the
 operator $e^{\sqrt{t}\Lambda_{1}}$ used in Theorem \ref{thm1.1}
is quantified by the operator $\Lambda_{1}$,
whose symbol is given by the $\ell^{1}$-norm $|\xi|_{1} =
|\xi_{1}|+|\xi_{2}|+|\xi_{3}|$, rather than the usual operator $\Lambda=\sqrt{-\Delta}$,
 whose
symbol is given by the $\ell^{2}$-norm $|\xi|= \sqrt{|\xi_{1}|^{2}+|\xi_{2}|^{2}+|\xi_{3}|^{2}}$,
as in Definition \ref{def1.1}.
 This approach enables us to avoid cumbersome
recursive estimation of higher order derivatives and intricate combinatorial arguments to get the desired
decay estimates of solutions, see \cite{BBT,AB,CHBW,PG}. We also emphasize that the operators $\Lambda_{1}$ and $\Lambda$
are equivalent as a Fourier multiplier.
\medskip

4.\  If we remove the exponential operator $e^{\sqrt{t}\Lambda_{1}}$,
 the result of Theorem \ref{thm1.1} is essentially
obtained in Liu-Zhang-Zhao \cite{LZZ} (see also Hao-Liu \cite{HL13}).
However, $e^{\sqrt{t}\Lambda_{1}}$ is of importance for
nonlinear estimates in the space $\widetilde{L}^{\rho}(0,T;e^{\sqrt{T}\Lambda_{1}}\dot{B}^{s}_{p,q})$,
 see Lemari\'{e}-Rieusset \cite{PG}. In fact, one of the main points in the proof of Theorem
 \ref{thm1.1} is how the product in the nonlinear terms of \eqref{eq1.1} is estimated in
 $\widetilde{L}^{1}(0,T^{*};e^{\sqrt{t}\Lambda_{1}}\dot{B}^{-1+\frac{3}{p}}_{p,1}(\mathbb{R}^{3}))$ or
$\widetilde{L}^{1}(0,T^{*};e^{\sqrt{t}\Lambda_{1}}\dot{B}^{\frac{3}{q}}_{q,1}(\mathbb{R}^{3}))$, although
the Plancherel
Theorem used in Foias-Temam \cite{FT} is no longer available,  we can resort
to the approach taken in \cite{PG} and \cite{BBT} to find out the nice boundedness
property of the following bilinear operator (cf. \eqref{eq3.12} below)
\begin{align*}
B_{t}(u,v)=&e^{\sqrt{t}\Lambda_{1}} (e^{-\sqrt{t}\Lambda_{1}} u e^{-\sqrt{t}\Lambda_{1}}v )
=\frac{1}{(2\pi)^{3}}\int_{\mathbb{R}^{3}}\int_{\mathbb{R}^{3}}e^{ix\cdot\xi} e^{\sqrt{t}(|\xi|_{1}-|\xi-\eta|_{1}-|\eta|_{1})}\widehat{u}(\xi-\eta)\widehat{v}(\eta)\text{d}\eta\text{d}\xi.
\end{align*}
Based on the desired property of $B_{t}(u,v)$, the estimates of the nonlinear terms can be done via
the  Fourier localization approach (see Lemmas \ref{lem3.3}, \ref{lem3.5} below).
\end{remark}


Notice that an important property of Gevrey regular functions is that estimates on higher-order derivatives follow immediately,
using this fact, it follows from Theorem \ref{thm1.1} and Stirling's approximation that the solution of
 system \eqref{eq1.1} with initial data $(u_{0},d_{0}-\overline{d}_{0})$ belonging
 to $\dot{B}^{\frac{3}{p}-1}_{p,1}(\mathbb{R}^{3})\times \dot{B}^{\frac{3}{q}}_{q,1}(\mathbb{R}^{3})$
 automatically satisfies certain higher-order temporal decay estimates.

\begin{corollary}\label{cor1.4}
Let $1<p, q<\infty$ satisfy \eqref{eq1.3},  $\overline{d}_{0}\in \mathbb{S}^{2}$ be a constant vector,
$u_{0}\in\dot{B}^{\frac{3}{p}-1}_{p,1}(\mathbb{R}^{3})$ with
$\operatorname{div} u_{0}=0$,
$d_{0}-\overline{d}_{0}\in\dot{B}^{\frac{3}{q}}_{q,1}(\mathbb{R}^{3})$,
and $M_{0}:=\|u_{0}\|_{\dot{B}^{\frac{3}{p}-1}_{p,1}}+\|d_{0}-\overline{d}_{0}\|_{\dot{B}^{\frac{3}{q}}_{q,1}}$.
 Let $(u,d)$  be the solution obtained in Theorem \ref{thm1.1}.
 Then we have for all $m\geq 0$ and $0<t<T^{*}$
\begin{align*}
\|\Lambda^{m}u(t)\|_{\widetilde{L}^{\infty}_{T^{*}}(\dot{B}^{\frac{3}{p}-1}_{p,1})
\cap \widetilde{L}^{1}_{T^{*}}(\dot{B}^{\frac{3}{p}+1}_{p,1})}+
\|\Lambda^{m}(d(t)-\overline{d}_{0})\|_{\widetilde{L}^{\infty}_{T^{*}}(\dot{B}^{\frac{3}{q}}_{q,1})
\cap \widetilde{L}^{1}_{T^{*}}(\dot{B}^{\frac{3}{q}+2}_{q,1})}\leq C t^{-\frac{m}{2}} M_{0},
\end{align*}
 where $C>0$ is an absolute constant. 
\end{corollary}

\textbf{The organization of the paper}. In Section 2, we collect some preliminary
materials, including the Littlewood-Paley decomposition, the definition of Besov spaces 
and some useful properties.
Section 3 is devoted to proving  Theorem \ref{thm1.1}.
In the last Section, we shall give the proof of Corollary \ref{cor1.4}.
\medskip

 Let us complete this section by describing the notation we shall use in this paper:
\medskip 
 
\textbf{Notation}. We denote by $C$ a generic positive
constant, which may  vary at different places. The notation $A\lesssim B$, we mean that there is an uniform constant $C$, which may be different on each lines, such that $A\leq C B$. Let $X$ be a Banach space. For $1\leq p\leq \infty$, the notation
$L^{p}_{T}(X)$ (or $L^{p}(0,T;X)$) stands for the set of measurable
functions on $(0,T)$ with values in $X$, such that $t\rightarrow
\|f(t)\|_{X}$ belongs to $L^{p}(0,T)$. For $x=(x_{1},x_{2},x_{3})\in \mathbb{R}^{3}$,
 we denote $|x|_{p}=\left(|x_{1}|^{p}+|x_{2}|^{p}+|x_{3}|^{p}\right)^{\frac{1}{p}}$ and $|x|=|x|_{2}$.

\section{Preliminaries}

In this section, we shall give the definition of the Besov spaces
and some useful lemmas. In order to define Besov spaces, we first
introduce the Littlewood--Paley decomposition theory. Let
$\mathcal{S}(\mathbb{R}^{3})$ be the Schwartz class of rapidly
decreasing functions, given $f\in \mathcal{S}(\mathbb{R}^{3})$, its
Fourier transform $\mathcal{F}f=\widehat{f}$ and its inverse Fourier
transform $\mathcal{F}^{-1}f=\check{f}$ are, respectively, defined
by
\begin{align*}
\widehat{f}(\xi):=\int_{\mathbb{R}^{3}}e^{-ix\cdot\xi}f(x)\text{d}x
\text{ and
 }\check{f}(x):=\frac{1}{(2\pi)^{3}}\int_{\mathbb{R}^{3}}e^{ix\cdot\xi}f(\xi)\text{d}\xi.
\end{align*}
 Let   $\chi, \varphi \in
\mathcal{S}(\mathbb{R}^{3})$ be two
nonnegative radial functions supported in $\mathfrak{B}=\{\xi\in
\mathbb{R}^{3}:|\xi|\leq \frac{4}{3}\}$ and
$\mathfrak{C}=\{\xi\in\mathbb{R}^{3}:\frac{3}{4}\leq |\xi|\leq
\frac{8}{3}\}$ respectively, such that
\begin{align*}
\sum_{j\in\mathbb{Z}}\varphi(2^{-j}\xi)=1 \quad\text{ for any
}\xi\in\mathbb{R}^{3}\backslash\{0\},\quad\text{ and
 }\quad \chi(\xi)+\sum_{j\geq 0}\varphi(2^{-j}\xi)=1\quad\text{ for any }
\xi\in\mathbb{R}^{3}.
\end{align*}
For $j\in \mathbb{Z}$, the homogeneous Littlewood--Paley projection
operators ${S}_{j}$ and ${\Delta}_{j}$ are, respectively,
defined as
\begin{align*}
{S}_{j}f=\chi(2^{-j}D)f=2^{3j}\int_{\mathbb{R}^{3}}\widetilde{h}(2^{j}y)f(x-y)\text{d}y,
\text{ where } \widetilde{h}=\mathcal{F}^{-1}\chi,
\end{align*}
and
\begin{align*}
{\Delta}_{j} f=\varphi(2^{-j}D)
f=2^{3j}\int_{\mathbb{R}^{3}}h(2^{-j}y)f(x-y)\text{d}y, \text{ where
} h=\mathcal{F}^{-1}\varphi.
\end{align*}
Informally, ${\Delta}_{j}$ is a  frequency projection to the
annulus $\{|\xi|\sim 2^{j}\}$, while ${S}_{j}$ is a frequency
projection to the ball $\{|\xi|\lesssim 2^{j} \}$. One can easily
verify that ${\Delta}_{j}{\Delta}_{k}f=0$ if $|j-k|\geq 2$,
and $\Delta_{j}(S_{k-1}f\Delta_{k}f)=0$ if $|j-k|\geq 5$. Let
\begin{align*}
\mathcal{S}_{h}:= \{\phi \in \mathcal{S}(\mathbb{R}^{3}),
\int_{\mathbb{R}^{3}}\phi(x)x^{\gamma}\text{d}x=0,
|\gamma|=0,1,2,\cdots \}.
\end{align*}
Then its dual is given by
\begin{align*}
\mathcal{S}'_{h}=\mathcal{S}'/\mathcal{S}^{\bot}_{h}=\mathcal{S}'/\mathcal{P},
\end{align*}
where $\mathcal{P}$ is the space of polynomial.

We now recall  the definitions of the stationary/time-dependent homogeneous Besov spaces from
\cite{BCD,PG}.

\begin{definition}\label{def2.1}
Let $s\in\mathbb{R}$,
$1\leq p,r\leq \infty$, we set
\begin{equation*}
\|u\|_{\dot{B}^{s}_{p,r}} \triangleq \left\{
\begin{array}{l}
\left(\sum_{j\in\mathbb{Z}}
2^{jsr}\|\Delta_{j}u\|_{L^{p}}^{r}\right)^{\frac{1}{r}}
\quad\text{ for } 1\leq r< \infty, \\
\sup_{j\in\mathbb{Z}}
2^{js}\|\Delta_{j}u\|_{L^{p}}\quad\quad\quad\text{ for }r=\infty.
\end{array}
\right.
\end{equation*}

\item$\centerdot$ For $s<\frac{3}{p}$ (or $s=\frac{3}{p}$ if $r=1$), we define
$\dot{B}^{s}_{p,r}(\mathbb{R}^{3})\triangleq
\{u\in\mathcal{S}'_{h}(\mathbb{R}^{3})\ | \ \|u\|_{\dot{B}^{s}_{p,r}}<\infty\}$.

\item$\centerdot$  If $k\in\mathbb{N}$ and $\frac{3}{p}+k\leq s<\frac{3}{p}+k+1$ (or $s=\frac{3}{p}+k+1$ if
$r=1$), then $\dot{B}^{s}_{p,r}(\mathbb{R}^{3})$ is defined as the subset
of distributions $u\in\mathcal{S}'(\mathbb{R}^{3})$ such that
$\partial^{\beta}u\in \dot{B}^{s-k}_{p,r}(\mathbb{R}^{3})$ whenever
$|\beta|=k$.

\end{definition}

\begin{definition}\label{def2.2}
 Let $ 0<T\leq \infty$,  $s\in\mathbb{R}$ and
$1\leq \rho, p,r\leq \infty$. We define the mixed time space
 $\widetilde{L}^{\rho}(0,T;\dot{B}^{s}_{p,r}(\mathbb{R}^{3}))$
as the completion of $\mathcal{C}([0,T],\mathcal{S}(\mathbb{R}^{3}))$ by the norm
\begin{equation*}
\|u\|_{\widetilde{L}^{\rho}_{T}(\dot{B}^{s}_{p,r})} \triangleq
\left(
\sum_{j\in\mathbb{Z}} 2^{jsr}\left(\int_{0}^{T}
\|\Delta_{j} u(\cdot,t)\|_{L^{p}}^{\rho}\text{d}t\right)^{\frac{1}{\rho}}\right)^{\frac{1}{r}}<\infty
\end{equation*}
with the usual change if $\rho=\infty$ or $r=\infty$. For simplicity, we will of some use
$\widetilde{L}^{\rho}_{T}(\dot{B}^{s}_{p,r}(\mathbb{R}^{3}))$
instead of $\widetilde{L}^{\rho}(0,T;\dot{B}^{s}_{p,r}(\mathbb{R}^{3}))$
later on.

\end{definition}

Let us recall the following basic facts on the Littlewood-Paley
theory from \cite{BCD,PG}, which will be used in the subsequent
sections.

\begin{lemma}\label{lem2.3}
Let $\mathcal{B}$ be a ball and $\mathcal{C}$  a ring of
$\mathbb{R}^{3}$. A constant $C$ exists so that for any positive
real number $\Lambda_{1}$, any nonnegative integer $k$  and any couple
of real numbers $(p,q)$ with $q\geq p\geq 1$, there hold
\begin{align*}
&\operatorname{supp} \widehat{u}\subset \Lambda_{1} \mathcal{B}\quad
\Rightarrow \quad \sup_{|\alpha|=k}\|\partial^{\alpha}
u\|_{L^{q}}\leq
C^{k+1}\Lambda_{1}^{k+3(\frac{1}{p}-\frac{1}{q})}\|u\|_{L^{p}},\\
&\operatorname{supp} \widehat{u}\subset \Lambda_{1} \mathcal{C}\quad
\Rightarrow \quad
C^{-1-k}\Lambda_{1}^{k}\|u\|_{L^{p}}\leq\sup_{|\alpha|=k}\|\partial^{\alpha}
u\|_{L^{p}}\leq
C^{k+1}\Lambda_{1}^{k}\|u\|_{L^{p}}.  
\end{align*}
\end{lemma}

We also notice that  functions in the homogeneous Besov space
$\dot{B}^{s}_{p,r}(\mathbb{R}^{3})$ have the following properties
(see \cite{BCD,PG}).

\begin{lemma}\label{lem2.4}
For $s\in\mathbb{R}$, $1\leq p,r\leq \infty$, then there hold\\
(1).\ The set $C_{0}^{\infty}(\mathbb{R}^{3})$ is dense in $\dot{B}^{s}_{p,r}(\mathbb{R}^{3})$ if $|s|\leq \frac{3}{p}$, $1\leq p,r<\infty$;\\
(2).\ For $1\leq p_{1}\leq p_{2}\leq \infty$ and $1\leq r_{1}\leq r_{2}\leq \infty$, we have the continuous imbedding  $\dot{B}^{s}_{p_{1},r_{1}}(\mathbb{R}^{3})\hookrightarrow \dot{B}^{s-3(\frac{1}{p_{1}}-\frac{1}{p_{2}})}_{p_{2},r_{2}}(\mathbb{R}^{3})$;\\
(3).\ $\|u\|_{\dot{B}^{s}_{p,r}}\approx \|\nabla u\|_{\dot{B}^{s-1}_{p,r}}$;\\
(4).\ Let $\Lambda=\sqrt{-\Delta}$ and $\sigma\in\mathbb{R}$, then the operator $\Lambda_{1}^{\sigma}$ is an isomorphism from $\dot{B}^{s}_{p,r}(\mathbb{R}^{3})$
to $\dot{B}^{s-\sigma}_{p,r}(\mathbb{R}^{3})$;\\
(5).\ For $s>0$, $\dot{B}^{s}_{p,r}(\mathbb{R}^{3})\cap L^{\infty}(\mathbb{R}^{3})$ is an algebra.  Moreover, $\dot{B}^{\frac{3}{p}}_{p,1}(\mathbb{R}^{3})$ is an algebra since there holds $\dot{B}^{\frac{3}{p}}_{p,1}(\mathbb{R}^{3})\hookrightarrow \dot{B}^{0}_{\infty,1}(\mathbb{R}^{3})\hookrightarrow L^{\infty}(\mathbb{R}^{3})$;\\
(6).\ For  $s_{1},s_{2}\in \mathbb{R}$ such that $s_{1}<s_{2}$ and
$\theta\in (0,1)$, we have the following interpolation inequalities:
\begin{align*}
&
\|u\|_{\dot{B}^{s_{1}\theta+s_{2}(1-\theta)}_{p,r}}\leq C\|u\|_{\dot{B}^{s_{1}}_{p,r}}^{\theta}\|u\|_{\dot{B}^{s_{2}}_{p,r}}^{1-\theta},\\
&
\|u\|_{\dot{B}^{s_{1}\theta+s_{2}(1-\theta)}_{p,1}}\leq \frac{C}{s_{2}-s_{1}}\left(\frac{1}{\theta}+\frac{1}{1-\theta}\right)\|u\|_{\dot{B}^{s_{1}}_{p,\infty}}^{\theta}
\|u\|_{\dot{B}^{s_{2}}_{p,\infty}}^{1-\theta},
\end{align*}
where $C$ is a positive constant.

\end{lemma}

Let us complete this section by describing the Bony's paraproduct decomposition.
Let $u$ and $v$ be two temperate distributions, the paraproduct
between $u$ and $v$ is defined by
\begin{align*}
T_{u}v \triangleq\sum_{j\in\mathbb{Z}}S_{j-1}u\Delta_{j}v.
\end{align*}
The remainder of the paraproduct $R(u,v)$ is defined by
\begin{align*}
R(u,v)\triangleq\sum_{|j-j'|\leq
1}\Delta_{j}u\Delta_{j'}v=\sum_{j\in\mathbb{Z}}\Delta_{j}u
\widetilde{\Delta}_{j} v, \quad \widetilde{\Delta}_{j} v\triangleq
\sum_{|j'-j|\leq 1} \Delta_{j}v.
\end{align*}
Then, we have the following Bony's decomposition:
\begin{align}\label{eq2.1}
uv=T_{u}v+T_{v}u+R(u,v).
\end{align}

\section{The proof of Theorem \ref{thm1.1}}

In this section, we shall give the proof of Theorem \ref{thm1.1}. In
order to do it, followed by some ideas as in \cite{BBT,PG,CHBW}, let us first establish
 the following two lemmas for
the basic heat equation
\begin{align}
\label{eq3.1}
\begin{cases}
\partial_{t}u-\Delta u =f  &\text{ in }\mathbb{R}^{3}\times (0,\infty),\\
 u|_{t=0}=u_{0}(x)    &\text{ in }\mathbb{R}^{3}.\\
\end{cases}
\end{align}
From the two Lemmas \ref{lem3.1} and \ref{lem3.2} below, we notice that
if we remove the operator $e^{\sqrt{t}\Lambda_{1}}$, the results were essentially obtained in
\cite{BCD} (see also \cite{HL13}). However, $e^{\sqrt{t}\Lambda_{1}}$ is of importance for the nonlinear
estimates, which are needed for the proof of our main result.

Now, by the Duhamel principle, we can express a solution $u$ of
\eqref{eq3.1} in the integral form:
\begin{align}\label{eq3.2}
u(t)= e^{ t\Delta} u_{0}+\int_{0}^{t} e^{(t-s)\Delta}
f(s)\text{d}s.
\end{align}

\begin{lemma}\label{lem3.1}
 For $u_{0}\in \dot{B}^{s}_{p,r}(\mathbb{R}^{3})$ with
$s\in\mathbb{R}$,
$1< p<\infty$ and $1\leq r\leq \infty$, we have\\
\ (1).\ Assume that $\|u_{0}\|_{\dot{B}^{s}_{p,r}}\leq C_{0}$,
$1\leq \rho<\infty$. For  any small $\varepsilon_{0}>0$, there
exists $T_{0}>0$, such that the following estimate holds
\begin{align}\label{eq3.3}
\|e^{ t\Delta+\sqrt{t}\Lambda_{1}}
u_{0}\|_{\widetilde{L}^{\rho}_{T}(\dot{B}^{s+\frac{2}{\rho}}_{p,r})}\leq
\varepsilon_{0} \text{ for any } 0<T\leq T_{0}.
\end{align}
\ (2).\ Assume that  $1\leq \rho\leq \infty$. For any small
$\varepsilon_{0}>0$, there exists $\eta_{0}$, such that if
$\|u_{0}\|_{\dot{B}^{s}_{p,r}}\leq \eta_{0}$, the following
estimate holds
\begin{align}\label{eq3.4}
\|e^{ t\Delta+\sqrt{t}\Lambda_{1}}u_{0}\|_{\widetilde{L}^{\rho}(\mathbb{R}_{+};
\dot{B}^{s+\frac{2}{\rho}}_{p,r})}\leq \varepsilon_{0}.
\end{align}
\end{lemma}

\begin{proof}
Before going to the proof, let us first recall that from Lemma \textbf{1} of Bae, Biswas and Tadmor \cite{BBT}
(see also Chapter 24, Lemari\'{e}-Rieusset \cite{PG}), the operator 
$e^{\frac{1}{2}t\Delta+\sqrt{t}\Lambda_{1}}$ is a Fourier multiplier
 which maps boundedly $L^{p}\rightarrow L^{p}$, $1<p<\infty$, and its o
 perator norm is uniformly bounded with respect to $t\geq 0$. Then, one easily sees that
\begin{align}\label{eq3.5}
\|\Delta_{j} e^{ t\Delta+\sqrt{t}\Lambda_{1}}
u_{0}\|_{L^{p}}\leq& \|\mathcal{F}^{-1} e^{-\frac{1}{2}t|\xi|^{2}+\sqrt{t}|\xi|_{1}}\|_{L^{1}} \|\Delta_{j}e^{ \frac{1}{2}t\Delta}
u_{0}\|_{L^{p}}\lesssim  \|\Delta_{j} e^{ \frac{1}{2}t\Delta}
u_{0}\|_{L^{p}}\nonumber\\
\lesssim & e^{-\kappa \frac{1}{2} 2^{2j}t} \|\Delta_{j} u_{0}\|_{L^{p}}.
\end{align}
Hence,
\begin{align*}
\|\Delta_{j} e^{ t\Delta+\sqrt{t}\Lambda_{1}}
u_{0}\|_{L^{\rho}_{T}L^{p}}
\lesssim & \left(\frac{1-e^{-\kappa \rho \frac{1}{2} 2^{2j}T}}{\kappa \rho \frac{1}{2}2^{2j}}\right)^{\frac{1}{\rho}} \|\Delta_{j} u_{0}\|_{L^{p}}.
\end{align*}
Noticing that we can choose a sufficiently large $J\in \mathbb{N}$ such that
\begin{align*}
\sum_{|j|>J}\left(\frac{1-e^{-\kappa \rho \frac{1}{2} 2^{2j}T}}{\kappa \rho \frac{1}{2}}\right)^{\frac{r}{\rho}} \left(2^{js}\|\Delta_{j} u_{0}\|_{L^{p}}\right)^{r}\leq \frac{\varepsilon_{0}}{2}.
\end{align*}
On the other hand, one can choose $T_{0}$ small enough such that for $0\leq T\leq T_{0}$,
\begin{align*}
\sum_{|j|\leq J}\left(\frac{1-e^{-\kappa \rho \frac{1}{2} 2^{2j}T}}{\kappa \rho \frac{1}{2}}\right)^{\frac{r}{\rho}}\left(2^{js} \|\Delta_{j} u_{0}\|_{L^{p}}\right)^{r}\leq \frac{\varepsilon_{0}}{2}.
\end{align*}
Combining the above two estimates together, we obtain \eqref{eq3.3}. To prove \eqref{eq3.4},
 multiplying both side of  \eqref{eq3.5} with  $2^{js}$,
  and then taking sequence $\ell^{r}$-norm in the resulting inequalities, it is easy to see that
\begin{align*}
\| e^{ t\Delta+\sqrt{t}\Lambda_{1}}
u_{0}\|_{\widetilde{L}^{\rho}_{T}(\dot{B}^{s+\frac{2}{\rho}}_{p,r})} \lesssim & \left\{\sum_{j\in\mathbb{Z}}\left(\frac{1-e^{-\kappa \rho \frac{1}{2} 2^{2j}T}}{\kappa \rho \frac{1}{2}}\right)^{\frac{r}{\rho}}\left(2^{js} \|\Delta_{j} u_{0}\|_{L^{p}}\right)^{r}\right\}^{\frac{1}{r}}\nonumber\\ \lesssim & \left\{\sum_{j\in\mathbb{Z}}\left(2^{js} \|\Delta_{j} u_{0}\|_{L^{p}}\right)^{r}\right\}^{\frac{1}{r}}\leq C \|u_{0}\|_{\dot{B}^{s}_{p,r}}.
\end{align*}
By taking $\eta_{0}\leq \frac{\varepsilon}{C}$, we obtain \eqref{eq3.4}.
 These complete the proof of Lemma \ref{lem3.1}.
\end{proof}

\begin{lemma}\label{lem3.2}
Let $T>0$, $s\in \mathbb{R}$, $1<p<\infty$ and $1\leq  \rho,
r\leq \infty$. Assume that $u_{0}\in \dot{B}^{s}_{p,r}(\mathbb{R}^{3})$ and $f\in
\widetilde{L}^{\rho}\!(0,T;e^{\!\sqrt{t}\Lambda_{1}}\dot{B}^{s-\!2+\!\frac{2}{\rho}}_{p,r}\!(\mathbb{R}^{3}))$. Then
\eqref{eq3.2} has a unique solution $u\!\in\!
\widetilde{L}^{\infty}(0,T; e^{\!\sqrt{t}\Lambda_{1}}\dot{B}^{s}_{p,r}\!(\mathbb{R}^{3}))
\cap \widetilde{L}^{\rho}(0,T; e^{\!\sqrt{t}\Lambda_{1}}\dot{B}^{s+\frac{2}{\rho}}_{p,r}\!(\mathbb{R}^{3}))$, and there exists a
constant $C$ such that for all $\rho_{1}\in [\rho,+\infty]$ such
that
\begin{align}\label{eq3.6}
\|u\|_{\widetilde{L}^{\rho_{1}}_{T}( e^{\sqrt{t}\Lambda_{1}}\dot{B}^{s+\frac{2}{\rho_{1}}}_{p,r})}
\leq C\left(\|u_{0}\|_{\dot{B}^{s}_{p,r}}+
\|f\|_{\widetilde{L}^{\rho}_{T}(e^{\sqrt{t}\Lambda_{1}}\dot{B}^{s-2+\frac{2}{\rho}}_{p,r})}\right).
\end{align}
\end{lemma}

\begin{proof}
In what follows, we shall show
\begin{align}\label{eq3.7}
\left\|\int_{0}^{t} e^{(t-s)\Delta}
f(s)\text{d}s\right\|_{\widetilde{L}^{\varrho_{1}}_{T}( e^{\sqrt{t}\Lambda_{1}}
\dot{B}^{s+\frac{2}{\varrho_{1}}}_{p,r})}
\lesssim \|f\|_{\widetilde{L}^{\rho}_{T}( e^{\sqrt{t}\Lambda_{1}}\dot{B}^{s-2+\frac{2}{\rho}}_{p,r})},
\end{align}
from which together with Lemma \ref{lem3.1}, it is  easy to  obtain \eqref{eq3.6}.
 Noticing that from Lemma \textbf{2} of Bae, Biswas and Tadmor \cite{BBT}(see also Chapter 24,
  Lemari\'{e}-Rieusset \cite{PG}),  the operator
  $e^{-(\sqrt{t-\tau}+\sqrt{\tau}-\sqrt{t})\Lambda_{1}}$
  for $0\leq \tau\leq t$ is either the identity operator
  or is an $L^{1}$ kernel whose $L^{1}$-norm is bounded independent of $\tau,t$. Then we have
\begin{align*}
\left\|\Delta_{j}\left(e^{\sqrt{t}\Lambda_{1}}\int_{0}^{t} e^{(t-\tau)\Delta}
f(\tau)\text{d}\tau\right)\right\|_{L^{p}}\lesssim& \int_{0}^{t} \left\| \Delta_{j}\left(e^{-(\sqrt{t-\tau}+\sqrt{\tau}-\sqrt{t})\Lambda_{1}}e^{\frac{1}{2}(t-\tau)\Delta+\sqrt{t-\tau}\Lambda_{1}}
e^{\frac{1}{2}(t-\tau)\Delta}e^{\sqrt{\tau}\Lambda_{1}}f(\tau)\right)\right\|_{L^{p}}\text{d}\tau\nonumber\\
\lesssim & \int_{0}^{t} \left\| \Delta_{j}\left(e^{\frac{1}{2}(t-\tau)\Delta+\sqrt{t-\tau}\Lambda_{1}}
e^{\frac{1}{2}(t-\tau)\Delta}e^{\sqrt{\tau}\Lambda_{1}}f(\tau)\right)\right\|_{L^{p}}\text{d}\tau\nonumber\\
\lesssim & \int_{0}^{t} \left\| \Delta_{j}\left(
e^{\frac{1}{2}(t-\tau)\Delta}e^{\sqrt{\tau}\Lambda_{1}}f(\tau)\right)\right\|_{L^{p}}\text{d}\tau\nonumber\\
\lesssim & \int_{0}^{t} e^{\kappa \frac{1}{2}(t-\tau) 2^{2j}}\left\| \Delta_{j}\left(
e^{\sqrt{\tau}\Lambda_{1}}f(\tau)\right)\right\|_{L^{p}}\text{d}\tau.
\end{align*}
By using Young's inequality, it follows that for $\rho\leq \rho_{1}$,
\begin{align*}
\left\|\Delta_{j}\left(e^{\sqrt{t}\Lambda_{1}}\int_{0}^{t} e^{(t-\tau)\Delta}
f(\tau)\text{d}\tau\right)\right\|_{L^{\rho_{1}}_{T}(L^{p})}
\lesssim 2^{-2j(1+\frac{1}{\rho_{1}}-\frac{1}{\rho})}
\left\|\Delta_{j} \left(e^{\sqrt{\tau}\Lambda_{1}} f(\tau)\right)\right\|_{L^{\rho}_{T}(L^{p})}.
\end{align*}
Multiplying both side of  the above inequality with  $2^{j(s+\frac{2}{\varrho_{1}})}$,
and then taking sequence $\ell^{r}$-norm in the resulting inequality, we obtain \eqref{eq3.7}.
This completes the proof of Lemma \ref{lem3.2}.
\end{proof}
\medskip

In what follows, we shall give
 two lemmas which play a crucial role in proof of our main result.

\begin{lemma}\label{lem3.3}
Let $0<T\leq \infty$ and $1<p,q<\infty$ such that
\begin{align}\label{eq3.8}
-\inf\{\frac{1}{3},\frac{1}{2p}\}\leq \frac{1}{q}-\frac{1}{p}.
\end{align}
Then, it holds that
\begin{align}\label{eq3.9}
\|fg\|_{\widetilde{L}^{1}_{T}( e^{\sqrt{t}\Lambda_{1}} \dot{B}^{\frac{3}{p}}_{p,1})}\lesssim
\|f\|_{\widetilde{L}^{\infty}_{T}( e^{\sqrt{t}\Lambda_{1}} \dot{B}^{\frac{3}{q}-1}_{q,1})}
\|g\|_{\widetilde{L}^{1}_{T}( e^{\sqrt{t}\Lambda_{1}} \dot{B}^{\frac{3}{q}+1}_{q,1})}
+\|f\|_{\widetilde{L}^{1}_{T}( e^{\sqrt{t}\Lambda_{1}} \dot{B}^{\frac{3}{q}+1}_{q,1})}
\|g\|_{\widetilde{L}^{\infty}_{T}( e^{\sqrt{t}\Lambda_{1}} \dot{B}^{\frac{3}{q}-1}_{q,1})},
\end{align}
for all $f,g\in \widetilde{L}^{\infty}_{T}( e^{\sqrt{t}\Lambda_{1}} \dot{B}^{\frac{3}{q}-1}_{q,1}(\mathbb{R}^{3}))\cap
\widetilde{L}^{1}_{T}( e^{\sqrt{t}\Lambda_{1}} \dot{B}^{\frac{3}{q}+1}_{q,1}(\mathbb{R}^{3}))$.
\end{lemma}

\begin{proof}
Set $F=e^{\sqrt{t}\Lambda_{1}}f$ and $G=e^{\sqrt{t}\Lambda_{1}}g$.
Then, by using Bony's paraproduct decomposition stated in \eqref{eq2.1} and the definitions of the  operators $\Delta_{j}$ and $S_{j}$, we get
\begin{align}\label{eq3.10}
e^{\sqrt{t}\Lambda_{1}}\Delta_{j}(fg)=& e^{\sqrt{t}\Lambda_{1}} \Delta_{j} (e^{-\sqrt{t}\Lambda_{1}}Fe^{-\sqrt{t}\Lambda_{1}}G)
\nonumber\\
=& \sum_{|j'-j|\leq 5}\sum_{j''\leq j'-2}\Delta_{j}( e^{-\sqrt{t}\Lambda_{1}} (\Delta_{j''}(e^{-\sqrt{t}\Lambda_{1}}F )\Delta_{j'}(e^{-\sqrt{t}\Lambda_{1}}G ) ))\nonumber\\
&+\sum_{|j'-j|\leq 5}\sum_{j''\leq j'-2}\Delta_{j}( e^{-\sqrt{t}\Lambda_{1}} (\Delta_{j''}(e^{-\sqrt{t}\Lambda_{1}}G )\Delta_{j'}(e^{-\sqrt{t}\Lambda_{1}}F ) ))\nonumber\\
&+\sum_{j'\geq j-N_{0}}\Delta_{j}(e^{-\sqrt{t}\Lambda_{1}} (\Delta_{j''}(e^{-\sqrt{t}\Lambda_{1}}F )\widetilde{\Delta}_{j'}(e^{-\sqrt{t}\Lambda_{1}}G  ))\nonumber\\
:=& I_{1}+I_{2}+I_{3}.
\end{align}
In what follows, we use the idea as in Lemari\'{e}-Rieusset \cite{PG} (see Chapter 24) and Bae-Biswas-Tadmor \cite{BBT} to consider the following operator $B_{t}(u,v)$ of the form:
\begin{align*}
B_{t}(u,v)=&e^{\sqrt{t}\Lambda_{1}} (e^{-\sqrt{t}\Lambda_{1}} u e^{-\sqrt{t}\Lambda_{1}}v )
=\frac{1}{(2\pi)^{3}}\int_{\mathbb{R}^{3}}\int_{\mathbb{R}^{3}}e^{ix\cdot\xi} e^{\sqrt{t}(|\xi|_{1}-|\xi-\eta|_{1}-|\eta|_{1})}\widehat{u}(\xi-\eta)\widehat{v}(\eta)\text{d}\eta\text{d}\xi.
\end{align*}
We now split the domain of integration of the above integral into sub-domains, depending on the sign of $\xi_{i}$, of $\eta_{i}$ and of $\xi_{i}-\eta_{i}$. Indeed,  denote for $\lambda=(\lambda_{1},
\lambda_{2},\lambda_{3})$, $\mu=(\mu_{1},\mu_{2},\mu_{3})$, $\gamma=(\gamma_{1},\gamma_{2},\gamma_{3})\in \{-1,1\}^{3}$,
\begin{align*}
&D_{\lambda}=\{\eta: \lambda_{i}\eta_{i}\geq 0, \quad i=1,2,3\};\quad D_{\mu}=\{\xi-\eta: \mu_{i}(\xi_{i}-\eta_{i})\geq 0, \quad i=1,2,3\};\\
&D_{\gamma}=\{\xi: \gamma_{i}\xi_{i}\geq 0, \quad i=1,2,3\}.
\end{align*}
Let $\chi_{D}$ be the characteristic function on the domain $D$. Then we can rewrite $B_{t}(u,v)$ as
\begin{align*}
B_{t}(u,v)\! =\!\!\sum_{\lambda,\mu,\gamma\in \{-1,1\}^{3}}\frac{1}{(2\pi)^{3}}\!\!
\int_{\!\mathbb{R}^{3}}\!\int_{\mathbb{R}^{3}}\! e^{ix\cdot\xi} \chi_{D_{\gamma}}(\xi) 
e^{\sqrt{t}(|\xi|_{1}-|\xi-\eta|_{1}-|\eta|_{1}) } \chi_{D_{\mu}}
(\xi-\eta)\widehat{u}(\xi-\eta)\chi_{D_{\lambda}}(\eta)\widehat{v}(\eta)\text{d}\eta\text{d}\xi.
\end{align*}
By using this observation, we introduce the following operators acting on one variable:
\begin{align*}
K_{k}u=
\begin{cases}
\frac{1}{2\pi} \int_{0}^{\infty} e^{ix\xi}\widehat{u}(\xi)\text{d}\xi \quad \text{ if } k=1;\\
\frac{1}{2\pi} \int_{-\infty}^{0} e^{ix\xi}\widehat{u}(\xi)\text{d}\xi\quad \text{ if } k=-1
\end{cases}
\end{align*}
and
\begin{align*}
\mathcal{L}_{t,k_{1},k_{2}} u =
\begin{cases}
u\quad \quad\quad \quad\quad \quad& \text{ if } k_{1},k_{2}\in \{-1,1\} \text{ and }k_{1}k_{2}=1;\\
\frac{1}{2\pi} \int_{\mathbb{R}}e^{ix\xi} 2^{-2t|\xi|}\widehat{u}(\xi)\text{d}\xi \quad &\text{ if }
 k_{1},k_{2}\in \{-1,1\} \text{ and }k_{1}k_{2}=-1.
\end{cases}
\end{align*}
For $\alpha=(\alpha_{1},\alpha_{2},\alpha_{3})$, $\beta=(\beta_{1},\beta_{2},\beta_{3})\in \{-1,1\}^{3}$, denote the operator
\begin{align*}
Z_{t,\alpha,\beta}=K_{\beta_{1}} \mathcal{L}_{t,\alpha_{1},\beta_{1}}\otimes\cdots \otimes K_{\beta_{3}} \mathcal{L}_{t,\alpha_{3},\beta_{3}}.
\end{align*}
It is easy to see that the above tensor product means that the $j$-th operator in the tensor product
acts on the $j$-th variable of the function $u(x_{1},x_{2},x_{3})$.
Then, a tedious but elementary calculations   yields the following identity
\begin{align}\label{eq3.11}
B_{t}(u,v)=\sum_{(\lambda,\mu,\gamma)\in \{-1,1\}^{3}}K_{\lambda_{1}}\otimes K_{\lambda_{2}}\otimes
K_{\lambda_{3}}(Z_{t,\lambda
,\mu} u Z_{t,\mu,\gamma} v).
\end{align}
Notice that for $\eta\in D_{\lambda}$, $\xi-\eta\in D_{\mu}$ and
$\xi\in D_{\gamma}$, $e^{\sqrt{t}(|\xi|_{1}-|\xi-\eta|_{1}-|\eta|_{1})}$ must belongs to the following set
\begin{align*}
\mathfrak{M}:=\{1,e^{-2\sqrt{t}|\xi_{i}|}, e^{-2\sqrt{t}|\xi_{i}-\eta_{i}|}e^{-2\sqrt{t}|\eta_{i}|}\}, \quad i=1,2,3.
\end{align*}
On the other hand, it is easy to see that the operators $\chi_{D_{\lambda}}$,
$\chi_{D_{\mu}}$, $\chi_{D_{\gamma}}$, and every element of $\mathfrak{M}$
are the Fourier multiplier on $L^{p}(\mathbb{R}^{3})$, $1<p<\infty$.
Hence, the operators $K_{\lambda_{i}}$ $(i=1,2,3)$, $Z_{t,\lambda,\mu }$
and $Z_{t,\mu,\gamma}$ in \eqref{eq3.11} are linear combination of Fourier
multipliers on $L^{p}(\mathbb{R}^{3})$ (including Hilbert transform)
and the identity operator. Moreover, they are bounded linear operators
on $L^{p}(\mathbb{R}^{3})$, $1<p<\infty$, and the corresponding norm of
$Z_{t,\lambda,\mu }$ and $Z_{t,\mu,\gamma}$ is bounded independent of $t\geq 0$.
By using these facts, we taking the $L^{p}$-norm to \eqref{eq3.11}, it follows that
\begin{align}\label{eq3.12}
\|B_{t}(u,v)\|_{L^{p}}\lesssim \|Z_{t,\lambda,\mu } u Z_{t,\mu,\gamma}v \|_{L^{p}}.
\end{align}
Applying this argument, we can estimate the term $I_{1}$ for
$1<p\leq \frac{3}{2}$ and $p\leq q\leq \frac{3p}{3-p}$,
\begin{align*}
\|I_{1}\|_{L^{1}_{T}L^{p}}\lesssim&
\sum_{|j'-j|\leq 5}\sum_{j''\leq j'-2}
\|\Delta_{j}( e^{-\sqrt{t}\Lambda_{1}}
(\Delta_{j''}(e^{-\sqrt{t}\Lambda_{1}}F )\Delta_{j'}(e^{-\sqrt{t}\Lambda_{1}}G ) ))\|_{L^{1}_{T}L^{p}}\nonumber\\
\lesssim& \sum_{|j'-j|\leq 5}\sum_{j''\leq j'-2}
\|\Delta_{j}B_{t}( \Delta_{j''}F,\Delta_{j'}G ) \|_{L^{1}_{T}L^{p}}
\lesssim \sum_{|j'-j|\leq 5}\sum_{j''\leq j'-2}
\|Z_{t,\lambda,\mu}\Delta_{j''}F Z_{t,\mu,\gamma}\Delta_{j'}G  \|_{L^{1}_{T}L^{p}}\nonumber\\
\lesssim& \sum_{|j'-j|\leq 5}\sum_{j''\leq j'-2}
\|Z_{t,\lambda,\mu}\Delta_{j''}F \|_{L^{\infty}_{T}L^{\frac{pq}{q-p}}}
\|Z_{t,\mu,\gamma}\Delta_{j'}G  \|_{L^{1}_{T}L^{q}}\nonumber\\
\lesssim&  \sum_{|j'-j|\leq 5}\sum_{j''\leq j'-2} 2^{(\frac{3}{q}-\frac{3(q-p)}{pq})j''}
 \|\Delta_{j''} Z_{t,\lambda,\mu}F \|_{L^{\infty}_{T}L^{q}}
  \|\Delta_{j'}Z_{t,\mu,\gamma} G  \|_{L^{1}_{T}L^{q}}\nonumber\\
\lesssim&  \sum_{|j'-j|\leq 5}\sum_{j''\leq j'-2} 2^{(1-\frac{3(q-p)}{pq})j''}
2^{(\frac{3}{q}-1)j''} \|\Delta_{j''} F \|_{L^{\infty}_{T}L^{q}}
 \|\Delta_{j'} G  \|_{L^{1}_{T}L^{q}}\nonumber\\
\lesssim& \|F \|_{L^{\infty}_{T}(\dot{B}^{\frac{3}{q}-1}_{q,1})}
\sum_{|j'-j|\leq 5}  2^{-\frac{3}{p}j'}2^{(1+\frac{3}{q})j'}
\|\Delta_{j'} G  \|_{L^{1}_{T}L^{q}}\nonumber\\
\lesssim& \|F \|_{L^{\infty}_{T}(\dot{B}^{\frac{3}{q}-1}_{q,1})}
\sum_{|j'-j|\leq 5}  2^{-\frac{3}{p}j'} a_{j'}
\| G  \|_{L^{1}_{T}(\dot{B}^{\frac{3}{q}+1}_{q,1})}\nonumber\\
\lesssim&  2^{-\frac{3}{p}j} a_{j}
 \|F \|_{L^{\infty}_{T}(\dot{B}^{\frac{3}{q}-1}_{q,1})} \| G  \|_{L^{1}_{T}(\dot{B}^{\frac{3}{q}+1}_{q,1})},
\end{align*}
for $1<q\leq p\leq \frac{3}{2}$,
\begin{align*}
\|I_{1}\|_{L^{1}_{T}L^{p}}\lesssim&
 \sum_{|j'-j|\leq 5}\sum_{j''\leq j'-2}
\|Z_{t,\lambda,\mu}\Delta_{j''}F Z_{t,\mu,\gamma}\Delta_{j'}G  \|_{L^{1}_{T}L^{p}}\nonumber\\
\lesssim& \sum_{|j'-j|\leq 5}\sum_{j''\leq j'-2}
\|Z_{t,\lambda,\mu}\Delta_{j''}F \|_{L^{\infty}_{T}L^{\infty}}
\|Z_{t,\mu,\gamma}\Delta_{j'}G  \|_{L^{1}_{T}L^{p}}\nonumber\\
\lesssim&  \sum_{|j'-j|\leq 5}2^{(\frac{3}{q}-\frac{3}{p})j'}\sum_{j''\leq j'-2} 2^{\frac{3}{q}j''}
 \|\Delta_{j''} Z_{t,\lambda,\mu}F \|_{L^{\infty}_{T}L^{q}}
  \|\Delta_{j'}Z_{t,\mu,\gamma} G  \|_{L^{1}_{T}L^{q}}\nonumber\\
\lesssim&  \sum_{|j'-j|\leq 5}2^{(\frac{3}{q}-\frac{3}{p})j'}
\sum_{j''\leq j'-2} 2^{j''}2^{(\frac{3}{q}-1)j''}
 \|\Delta_{j''} F \|_{L^{\infty}_{T}L^{q}}
 \|\Delta_{j'} G  \|_{L^{1}_{T}L^{q}}\nonumber\\
\lesssim& \|F \|_{L^{\infty}_{T}(\dot{B}^{\frac{3}{q}-1}_{q,1})}
\sum_{|j'-j|\leq 5}  2^{-\frac{3}{p}j'}2^{(1+\frac{3}{q})j'}
\|\Delta_{j'} G  \|_{L^{1}_{T}L^{q}}\nonumber\\
\lesssim&  2^{-\frac{3}{p}j} a_{j}
 \|F \|_{L^{\infty}_{T}(\dot{B}^{\frac{3}{q}-1}_{q,1})} \| G  \|_{L^{1}_{T}(\dot{B}^{\frac{3}{q}+1}_{q,1})},
\end{align*}
and for $\frac{3}{2}\leq p<\infty$ and $1<q\leq 2p$,
\begin{align*}
\|I_{1}\|_{L^{1}_{T}L^{p}}\lesssim&
 \sum_{|j'-j|\leq 5}\sum_{j''\leq j'-2}
\|Z_{t,\lambda,\mu}\Delta_{j''}F Z_{t,\mu,\gamma}\Delta_{j'}G  \|_{L^{1}_{T}L^{p}}\nonumber\\
\lesssim& \sum_{|j'-j|\leq 5}\sum_{j''\leq j'-2}
\|Z_{t,\lambda,\mu}\Delta_{j''}F\|_{L^{\infty}_{T}L^{2p}}
\| Z_{t,\mu,\gamma}\Delta_{j'}G  \|_{L^{1}_{T}L^{2p}}\nonumber\\
\lesssim&  \sum_{|j'-j|\leq 5}2^{(\frac{3}{q}-\frac{3}{2p})j'}\sum_{j''\leq j'-2}
 2^{(\frac{3}{q}-\frac{3}{2p})j''}
 \|\Delta_{j''} Z_{t,\lambda,\mu}F \|_{L^{\infty}_{T}L^{q}}
  \|\Delta_{j'}Z_{t,\mu,\gamma} G  \|_{L^{1}_{T}L^{q}}\nonumber\\
\lesssim&  \sum_{|j'-j|\leq 5}2^{(\frac{3}{q}-\frac{3}{2p})j'}
\sum_{j''\leq j'-2} 2^{(1-\frac{3}{2p})j''}2^{(\frac{3}{q}-1)j''}
 \|\Delta_{j''} F \|_{L^{\infty}_{T}L^{q}}
 \|\Delta_{j'} G  \|_{L^{1}_{T}L^{q}}\nonumber\\
\lesssim& \|F \|_{L^{\infty}_{T}(\dot{B}^{\frac{3}{q}-1}_{q,1})}
\sum_{|j'-j|\leq 5}  2^{-\frac{3}{p}j'}2^{(1+\frac{3}{q})j'}
\|\Delta_{j'} G  \|_{L^{1}_{T}L^{q}}\nonumber\\
\lesssim&  2^{-\frac{3}{p}j} a_{j}
 \|F \|_{L^{\infty}_{T}(\dot{B}^{\frac{3}{q}-1}_{q,1})} \| G  \|_{L^{1}_{T}(\dot{B}^{\frac{3}{q}+1}_{q,1})},
\end{align*}
where we have used the fact that $Z_{t,\lambda,\nu}$ commutes with $\Delta_{j}$ and the boundedness
of $Z_{t,\lambda,\mu}$ on $L^{p}(\mathbb{R}^{3})$, $1<p<\infty$.
Here, $\{a_{j}\}_{j\in\mathbb{Z}}$ is the generic element of $\ell^{1}(\mathbb{Z})$ so that
$a_{j}\geq 0$ and $\sum_{j\in\mathbb{Z}}a_{j}=1$.  Similarly, for all $1<p,q<\infty$ and $-\inf\{\frac{1}{3},\frac{1}{2p}\}
\leq \frac{1}{q}-\frac{1}{p}$, we have
\begin{align*}
\|I_{2}\|_{L^{1}_{T}L^{p}}\lesssim&
 2^{-\frac{3}{p}j} a_{j}
 \|G \|_{L^{\infty}_{T}(\dot{B}^{\frac{3}{q}-1}_{q,1})} \| F  \|_{L^{1}_{T}(\dot{B}^{\frac{3}{q}+1}_{q,1})}.
\end{align*}
For the remainder term $I_{3}$, notice that when $1<p\leq \frac{3}{2}$,
we still have $q\leq \frac{3p}{3-p}\leq 2p$. Hence
\begin{align*}
\|I_{3}\|_{L^{1}_{T}L^{p}}\lesssim&
 \sum_{j'\geq j-N_{0}}\|\Delta_{j}(e^{-\sqrt{t}\Lambda_{1}} (\Delta_{j''}(e^{-\sqrt{t}\Lambda_{1}}F )
 \widetilde{\Delta}_{j'}(e^{-\sqrt{t}\Lambda_{1}}G  ))\|_{L^{1}_{T}L^{p}}\nonumber\\
 \lesssim&
 \sum_{j'\geq j-N_{0}}
 \|\Delta_{j}B_{t}( \Delta_{j'}F,\widetilde{\Delta}_{j'}G ) \|_{L^{1}_{T}L^{p}}
 \lesssim  \sum_{j'\geq j-N_{0}}
 \|Z_{t,\lambda,\mu}\Delta_{j'}F Z_{t,\mu,\gamma}\widetilde{\Delta}_{j'}G
  \|_{L^{1}_{T}L^{p}}\nonumber\\
 \lesssim&  \sum_{j'\geq j-N_{0}}\|Z_{t,\lambda,\mu}\Delta_{j'}F\|_{L^{1}_{T}L^{2p}}
 \|Z_{t,\mu,\gamma}\widetilde{\Delta}_{j'}G  \|_{L^{\infty}_{T}L^{2p}}
 \lesssim \sum_{j'\geq j-N_{0}}\|\Delta_{j'}F\|_{L^{1}_{T}L^{2p}}
 \|\widetilde{\Delta}_{j'}G  \|_{L^{\infty}_{T}L^{2p}}\nonumber\\
  \lesssim& \sum_{j'\geq j-N_{0}} 2^{2(\frac{3}{q}-\frac{3}{2p})j'}
  \|\Delta_{j'}F\|_{L^{1}_{T}L^{q}}
 \|\widetilde{\Delta}_{j'}G  \|_{L^{\infty}_{T}L^{q}}\nonumber\\
 \lesssim&  \sum_{j'\geq j-N_{0}} 2^{(1+\frac{3}{q}-\frac{3}{p})j'} \|\Delta_{j'}F\|_{L^{1}_{T}L^{q}}
 \|G  \|_{L^{\infty}_{T}(\dot{B}^{\frac{3}{q}-1}_{q,1})} \nonumber\\
 \lesssim& 2^{-\frac{3}{p}j} a_{j}\| F  \|_{L^{1}_{T}(\dot{B}^{\frac{3}{q}+1}_{q,1})}
 \|G \|_{L^{\infty}_{T}(\dot{B}^{\frac{3}{q}-1}_{q,1})} .
\end{align*}
  Combining the estimates $I_{1}$, $I_{2}$, $I_{3}$ and the equality \eqref{eq3.10} together,
   we immediately obtain \eqref{eq3.9}. This completes the proof of Lemma \ref{lem3.3}.
\end{proof}

\begin{remark}\label{rem3.4}
We notice that if we replace \eqref{eq3.8} by
$-\inf\{\frac{1}{3},\frac{1}{2p}\}<\frac{1}{q}-\frac{1}{p}$ in Lemma
\ref{lem3.3},  with suitable revision of the proofs presented above,  we have for some $\theta\in (0,1)$
\begin{align*}
&\|f g\|_{\widetilde{L}^{1}_{T}( e^{\sqrt{t}\Lambda_{1}}\dot{B}^{\frac{3}{p}-1}_{p,1})}
\!\lesssim\!
\|f\|_{\widetilde{L}^{\frac{1+\theta}{\theta}}_{T}( e^{\sqrt{t}\Lambda_{1}}\dot{B}^{\frac{3}{q}-\frac{1-\theta}{1+\theta}}_{q,1})}
\|g\|_{\widetilde{L}^{1+\theta}_{T}( e^{\sqrt{t}\Lambda_{1}}\dot{B}^{\frac{3}{q}+\frac{1+\theta}{1+\theta}}_{q,1})}
\!\!+\!
\|f\|_{\widetilde{L}^{1+\theta}_{T}( e^{\sqrt{t}\Lambda_{1}}\dot{B}^{\frac{3}{q}+\frac{1+\theta}{1+\theta}}_{q,1})}
\|g\|_{\widetilde{L}^{\frac{1+\theta}{\theta}}_{T}( e^{\sqrt{t}\Lambda_{1}}\dot{B}^{\frac{3}{q}-\frac{1-\theta}{1+\theta}}_{q,1})}.
\end{align*}
\end{remark}

\begin{lemma}\label{lem3.5}
Let $0<T\leq \infty$ and $1<p,q<\infty$ such that
\begin{align}\label{eq3.13}
\frac{1}{q}-\frac{1}{p}\leq \frac{1}{3}.
\end{align}
Then, it holds that
\begin{align}\label{eq3.14}
\|f g\|_{\widetilde{L}^{1}_{T}( e^{\sqrt{t}\Lambda_{1}} \dot{B}^{\frac{3}{q}}_{q,1})}\lesssim
\|f\|_{\widetilde{L}^{2}_{T}( e^{\sqrt{t}\Lambda_{1}} \dot{B}^{\frac{3}{p}}_{p,1})}
\|g\|_{\widetilde{L}^{2}_{T}( e^{\sqrt{t}\Lambda_{1}} \dot{B}^{\frac{3}{q}}_{q,1})}
+\|f\|_{\widetilde{L}^{1}_{T}( e^{\sqrt{t}\Lambda_{1}} \dot{B}^{\frac{3}{p}+1}_{p,1})}
\|g\|_{\widetilde{L}^{\infty}_{T}( e^{\sqrt{t}\Lambda_{1}} \dot{B}^{\frac{3}{q}-1}_{q,1})},
\end{align}
for $f\in \widetilde{L}^{\infty}_{T}( e^{\sqrt{t}\Lambda_{1}} \dot{B}^{\frac{3}{p}-1}_{p,1}(\mathbb{R}^{3}))\cap
\widetilde{L}^{1}_{T}( e^{\sqrt{t}\Lambda_{1}} \dot{B}^{\frac{3}{p}+1}_{p,1}(\mathbb{R}^{3}))$ and
 $g \in \widetilde{L}^{\infty}_{T}( e^{\sqrt{t}\Lambda_{1}} \dot{B}^{\frac{3}{q}-1}_{q,1}(\mathbb{R}^{3}))\cap
\widetilde{L}^{1}_{T}( e^{\sqrt{t}\Lambda_{1}} \dot{B}^{\frac{3}{q}+1}_{q,1}(\mathbb{R}^{3}))$.
\end{lemma}

\begin{proof}
Set $F=e^{\sqrt{t}\Lambda_{1}}f$ and $G=e^{\sqrt{t}\Lambda_{1}}g$. It is easy to see that the identity \eqref{eq3.8}
 still holds.
In what follows, we only need to estimate the terms $I_{1}$, $I_{2}$ and $I_{3}$ of \eqref{eq3.8}.
By using \eqref{eq3.12} and the H\"{o}lder's inequality, $I_{1}$ can be estimated as
\begin{align*}
\|I_{1}\|_{L^{1}_{T}L^{q}}\!
\lesssim&\! \sum_{|j'-j|\leq 5}\sum_{j''\leq j'-2}
\|\Delta_{j}B_{t}( \Delta_{j''}F,\Delta_{j'}G ) \|_{L^{1}_{T}L^{q}}
\lesssim \!\!\sum_{|j'-j|\leq 5}\sum_{j''\leq j'-2}
\|Z_{t,\lambda,\mu}\Delta_{j''}F Z_{t,\mu,\gamma}\Delta_{j'}G  \|_{L^{1}_{T}L^{q}}\nonumber\\
\lesssim& \sum_{|j'-j|\leq 5}\sum_{j''\leq j'-2}
\|Z_{t,\lambda,\mu}\Delta_{j''}F \|_{L^{2}_{T}L^{\infty}}
\|Z_{t,\mu,\gamma}\Delta_{j'}G  \|_{L^{2}_{T}L^{q}}\nonumber\\
\lesssim&  \sum_{|j'-j|\leq 5}\sum_{j''\leq j'-2} 2^{\frac{3}{p}j''}
 \|\Delta_{j''} Z_{t,\lambda,\mu}F \|_{L^{2}_{T}L^{p}}
  \|\Delta_{j'}Z_{t,\mu,\gamma} G  \|_{L^{2}_{T}L^{q}}\nonumber\\
\lesssim& \|F \|_{L^{2}_{T}(\dot{B}^{\frac{3}{p}}_{p,1})}
\sum_{|j'-j|\leq 5}  2^{-\frac{3}{q}j'}2^{\frac{3}{q}j'}
\|\Delta_{j'} G  \|_{L^{2}_{T}L^{q}}
\lesssim  2^{-\frac{3}{q}j} a_{j}
 \|F \|_{L^{2}_{T}(\dot{B}^{\frac{3}{p}}_{p,1})} \| G  \|_{L^{2}_{T}(\dot{B}^{\frac{3}{q}}_{q,1})}.
\end{align*}
For the term $I_{2}$, notice that there holds $\frac{1}{q}-\frac{1}{p}\leq \frac{1}{3}$, we have  for $q\leq p$,
\begin{align*}
\|I_{2}\|_{L^{1}_{T}L^{q}}\!
\lesssim&\! \sum_{|j'-j|\leq 5}\sum_{j''\leq j'-2}
\|\Delta_{j}B_{t}( \Delta_{j''}G,\Delta_{j'}F ) \|_{L^{1}_{T}L^{q}}
\lesssim \!\!\sum_{|j'-j|\leq 5}\sum_{j''\leq j'-2}
\|Z_{t,\lambda,\mu}\Delta_{j''}G Z_{t,\mu,\gamma}\Delta_{j'}F  \|_{L^{1}_{T}L^{q}}\nonumber\\
\lesssim& \sum_{|j'-j|\leq 5}\sum_{j''\leq j'-2}
\|Z_{t,\lambda,\mu}\Delta_{j''}G \|_{L^{\infty}_{T}L^{r}}
\|Z_{t,\mu,\gamma}\Delta_{j'}F  \|_{L^{1}_{T}L^{p}}
\quad (\text{with }\frac{1}{q}=\frac{1}{p}+\frac{1}{r})\nonumber\\
\lesssim&  \sum_{|j'-j|\leq 5}\sum_{j''\leq j'-2} 2^{(\frac{3}{q}-\frac{3}{r})j''}
 \|\Delta_{j''} Z_{t,\lambda,\mu}G \|_{L^{\infty}_{T}L^{q}}
  \|\Delta_{j'}Z_{t,\mu,\gamma} F  \|_{L^{1}_{T}L^{p}}\nonumber\\
\lesssim&  \sum_{|j'-j|\leq 5}\sum_{j''\leq j'-2} 2^{(1-\frac{3}{r})j''}
2^{(\frac{3}{q}-1)j''} \|\Delta_{j''} G \|_{L^{\infty}_{T}L^{q}}
 \|\Delta_{j'} F  \|_{L^{1}_{T}L^{p}}
 \quad (\frac{1}{q}\!-\!\frac{1}{p}\leq \!\frac{1}{3} \Rightarrow 1\!-\!\frac{3}{r}\!\geq\! 0)\nonumber\\
\lesssim& \|G \|_{L^{\infty}_{T}(\dot{B}^{\frac{3}{q}-1}_{q,1})}
\sum_{|j'-j|\leq 5}  2^{-\frac{3}{q}j'}2^{(1+\frac{3}{p})j'}
\|\Delta_{j'} F  \|_{L^{1}_{T}L^{p}}\quad (1-\frac{3}{r}=1-\frac{3}{q}+\frac{3}{p})\nonumber\\
\lesssim&  2^{-\frac{3}{q}j} a_{j}
 \|G \|_{L^{\infty}_{T}(\dot{B}^{\frac{3}{q}-1}_{q,1})} \| F  \|_{L^{1}_{T}(\dot{B}^{\frac{3}{p}+1}_{p,1})},
\end{align*}
and for $q\geq p$,
\begin{align*}
\|I_{2}\|_{L^{1}_{T}L^{q}}\!
\lesssim&\! 2^{(\frac{3}{p}-\frac{3}{q})j}\sum_{|j'-j|\leq 5}\sum_{j''\leq j'-2}
\|\Delta_{j}B_{t}( \Delta_{j''}G,\Delta_{j'}F ) \|_{L^{1}_{T}L^{p}}\nonumber\\
\lesssim &2^{(\frac{3}{p}-\frac{3}{q})j}\sum_{|j'-j|\leq 5}\sum_{j''\leq j'-2}
\|Z_{t,\lambda,\mu}\Delta_{j''}G Z_{t,\mu,\gamma}\Delta_{j'}F  \|_{L^{1}_{T}L^{p}}\nonumber\\
\lesssim& 2^{(\frac{3}{p}-\frac{3}{q})j}\sum_{|j'-j|\leq 5}\sum_{j''\leq j'-2}
\|Z_{t,\lambda,\mu}\Delta_{j''}G \|_{L^{\infty}_{T}L^{\infty}}
\|Z_{t,\mu,\gamma}\Delta_{j'}F  \|_{L^{1}_{T}L^{p}}\nonumber\\
\lesssim& 2^{(\frac{3}{p}-\frac{3}{q})j}
 \sum_{|j'-j|\leq 5}\sum_{j''\leq j'-2} 2^{\frac{3}{q}j''}
 \|\Delta_{j''} Z_{t,\lambda,\mu}G \|_{L^{\infty}_{T}L^{q}}
  \|\Delta_{j'}Z_{t,\mu,\gamma} F  \|_{L^{1}_{T}L^{p}}\nonumber\\
\lesssim& 2^{(\frac{3}{p}-\frac{3}{q})j} \sum_{|j'-j|\leq 5}\sum_{j''\leq j'-2} 2^{j''}
2^{(\frac{3}{q}-1)j''} \|\Delta_{j''} G \|_{L^{\infty}_{T}L^{q}}
 \|\Delta_{j'} F  \|_{L^{1}_{T}L^{p}}\nonumber\\
\lesssim& 2^{(\frac{3}{p}-\frac{3}{q})j} \|G \|_{L^{\infty}_{T}(\dot{B}^{\frac{3}{q}-1}_{q,1})}
\sum_{|j'-j|\leq 5}  2^{-\frac{3}{p}j'}2^{(1+\frac{3}{p})j'}
\|\Delta_{j'} F  \|_{L^{1}_{T}L^{p}}\nonumber\\
\lesssim&  2^{-\frac{3}{q}j} a_{j}
 \|G \|_{L^{\infty}_{T}(\dot{B}^{\frac{3}{q}-1}_{q,1})} \| F  \|_{L^{1}_{T}(\dot{B}^{\frac{3}{p}+1}_{p,1})}.
\end{align*}
For the remainder term $I_{3}$, if in the case when $\frac{1}{p}+\frac{1}{q}\leq 1$.
Let $\frac{1}{r}=\frac{1}{p}+\frac{1}{q}$, then it holds that
\begin{align*}
\|I_{3}\|_{L^{1}_{T}L^{q}}\lesssim& 2^{\frac{3}{p}j}
 \sum_{j'\geq j-N_{0}}\|\Delta_{j}(e^{-\sqrt{t}\Lambda_{1}} (\Delta_{j''}(e^{-\sqrt{t}\Lambda_{1}}F )
 \widetilde{\Delta}_{j'}(e^{-\sqrt{t}\Lambda_{1}}G  ))\|_{L^{1}_{T}L^{r}}\nonumber\\
 \lesssim& 2^{\frac{3}{p}j}
 \sum_{j'\geq j-N_{0}}
 \|\Delta_{j}B_{t}( \Delta_{j'}F,\widetilde{\Delta}_{j'}G ) \|_{L^{1}_{T}L^{r}}
 \lesssim  2^{\frac{3}{p}j} \sum_{j'\geq j-N_{0}}
 \|Z_{t,\lambda,\mu}\Delta_{j'}F Z_{t,\mu,\gamma}\widetilde{\Delta}_{j'}G
  \|_{L^{1}_{T}L^{r}}\nonumber\\
 \lesssim&  2^{\frac{3}{p}j} \!\!\sum_{j'\geq j-N_{0}}\!\!
 \|Z_{t,\lambda,\mu}\Delta_{j'}F\|_{L^{2}_{T}L^{p}}
 \|Z_{t,\mu,\gamma}\widetilde{\Delta}_{j'}G  \|_{L^{2}_{T}L^{q}}\!
 \lesssim 2^{\frac{3}{p}j}\!\! \sum_{j'\geq j-N_{0}}\!\!\|\Delta_{j'}F\|_{L^{2}_{T}L^{p}}
 \|\widetilde{\Delta}_{j'}G  \|_{L^{2}_{T}L^{q}}\nonumber\\
  \lesssim& 2^{\frac{3}{p}j} \sum_{j'\geq j-N_{0}} 2^{-(\frac{3}{q}+\frac{3}{p})j'} 2^{\frac{3}{p}j'}
  \|\Delta_{j'}F\|_{L^{2}_{T}L^{p}}
 \|G  \|_{L^{2}_{T}(\dot{B}^{\frac{3}{q}}_{q,1})}\nonumber\\
 \lesssim&  2^{-\frac{3}{q}j}a_{j}\| F  \|_{L^{2}_{T}(\dot{B}^{\frac{3}{p}}_{p,1})}
 \|G  \|_{L^{2}_{T}(\dot{B}^{\frac{3}{q}}_{q,1})}.
\end{align*}
If in the case when $\frac{1}{p}+\frac{1}{q}> 1=\frac{1}{r}+\frac{1}{q}$, it holds that
\begin{align*}
\|I_{3}\|_{L^{1}_{T}L^{q}}\lesssim& 2^{3(1-\frac{1}{q})j}
 \sum_{j'\geq j-N_{0}}\|\Delta_{j}(e^{-\sqrt{t}\Lambda_{1}} (\Delta_{j''}(e^{-\sqrt{t}\Lambda_{1}}F )
 \widetilde{\Delta}_{j'}(e^{-\sqrt{t}\Lambda_{1}}G  ))\|_{L^{1}_{T}L^{1}}\nonumber\\
 \lesssim& 2^{3(1-\frac{1}{q})j}
 \sum_{j'\geq j-N_{0}}
 \|\Delta_{j}B_{t}( \Delta_{j'}F,\widetilde{\Delta}_{j'}G ) \|_{L^{1}_{T}L^{1}}
 \lesssim  2^{3(1-\frac{1}{q})j} \sum_{j'\geq j-N_{0}}
 \|Z_{t,\lambda,\mu}\Delta_{j'}F Z_{t,\mu,\gamma}\widetilde{\Delta}_{j'}G
  \|_{L^{1}_{T}L^{1}}\nonumber\\
 \lesssim&  2^{3(1-\frac{1}{q})j} \!\!\sum_{j'\geq j-N_{0}}
 \!\!\|Z_{t,\lambda,\mu}\Delta_{j'}F\|_{L^{2}_{T}L^{r}}
 \|Z_{t,\mu,\gamma}\widetilde{\Delta}_{j'}G  \|_{L^{2}_{T}L^{q}}\!
 \lesssim 2^{3(1-\frac{1}{q})j}\!\! \sum_{j'\geq j-N_{0}}\!\!\|\Delta_{j'}F\|_{L^{2}_{T}L^{r}}
 \|\widetilde{\Delta}_{j'}G  \|_{L^{2}_{T}L^{q}}\nonumber\\
  \lesssim& 2^{3(1-\frac{1}{q})j} \sum_{j'\geq j-N_{0}}
  2^{-(\frac{3}{q}+\frac{3}{r})j'} 2^{\frac{3}{p}j'}
  \|\Delta_{j'}F\|_{L^{2}_{T}L^{p}}
 \|G  \|_{L^{2}_{T}(\dot{B}^{\frac{3}{q}}_{q,1})}
 \lesssim 2^{-\frac{3}{q}j}a_{j}\| F  \|_{L^{2}_{T}(\dot{B}^{\frac{3}{p}}_{p,1})}
 \|G  \|_{L^{2}_{T}(\dot{B}^{\frac{3}{q}}_{q,1})}.
\end{align*}
  Combining the estimates $I_{1}$, $I_{2}$, $I_{3}$ together, it follows that \eqref{eq3.14} is established.
 This completes the proof of Lemma \ref{lem3.5}.
\end{proof}

\begin{remark}\label{rem3.6}
Notice that it holds that
$\|\nabla g\|_{\widetilde{L}^{\rho}_{T}( e^{\sqrt{t}\Lambda_{1}} \dot{B}^{s-1}_{p,r})}\approx
\| g\|_{\widetilde{L}^{\rho}_{T}( e^{\sqrt{t}\Lambda_{1}} \dot{B}^{s}_{p,r})}$.
From Lemma \ref{lem3.5}, it is easy to see that for
$f\in \widetilde{L}^{\infty}_{T}( e^{\sqrt{t}\Lambda_{1}} \dot{B}^{\frac{3}{p}-1}_{p,1}(\mathbb{R}^{3}))\cap
\widetilde{L}^{1}_{T}( e^{\sqrt{t}\Lambda_{1}} \dot{B}^{\frac{3}{p}+1}_{p,1}(\mathbb{R}^{3}))$ and
 $g \in \widetilde{L}^{\infty}_{T}( e^{\sqrt{t}\Lambda_{1}} \dot{B}^{\frac{3}{q}}_{q,1}(\mathbb{R}^{3}))\cap
\widetilde{L}^{1}_{T}( e^{\sqrt{t}\Lambda_{1}} \dot{B}^{\frac{3}{q}+2}_{q,1}(\mathbb{R}^{3}))$,
we have for  $\frac{1}{q}-\frac{1}{p}\leq \frac{1}{3}$
\begin{align*}
\|f\cdot\nabla g\|_{\widetilde{L}^{1}_{T}( e^{\sqrt{t}\Lambda_{1}} \dot{B}^{\frac{3}{q}}_{q,1})}\lesssim
\|f\|_{\widetilde{L}^{2}_{T}( e^{\sqrt{t}\Lambda_{1}} \dot{B}^{\frac{3}{p}}_{p,1})}
\|g\|_{\widetilde{L}^{2}_{T}( e^{\sqrt{t}\Lambda_{1}} \dot{B}^{\frac{3}{q}+1}_{q,1})}
+\|f\|_{\widetilde{L}^{1}_{T}( e^{\sqrt{t}\Lambda_{1}} \dot{B}^{\frac{3}{p}+1}_{p,1})}
\|g\|_{\widetilde{L}^{\infty}_{T}( e^{\sqrt{t}\Lambda_{1}} \dot{B}^{\frac{3}{q}}_{q,1})}.
\end{align*}
Moreover, if we replace \eqref{eq3.13} by $\frac{1}{q}-\frac{1}{p}<\frac{1}{3}$
in Lemma \ref{lem3.5}, with little revision of the proofs,   we have for some $\theta\in (0,1)$
\begin{align*}
\|f\cdot\nabla
g\|_{\widetilde{L}^{1}_{T}( e^{\sqrt{t}\Lambda_{1}}\dot{B}^{\frac{3}{q}}_{q,1})}
\lesssim
\|f\|_{\widetilde{L}^{2}_{T}( e^{\sqrt{t}\Lambda_{1}}\dot{B}^{\frac{3}{p}}_{p,1})}
\|g\|_{\widetilde{L}^{2}_{T}( e^{\sqrt{t}\Lambda_{1}}\dot{B}^{\frac{3}{q}+1}_{q,1})}
+\|f\|_{\widetilde{L}^{1+\theta}_{T}( e^{\sqrt{t}\Lambda_{1}}\dot{B}^{\frac{3}{p}+\frac{2}{1+\theta}}_{p,1})}
\|g\|_{\widetilde{L}^{\frac{1+\theta}{\theta}}_{T}( e^{\sqrt{t}\Lambda_{1}}
\dot{B}^{\frac{3}{q}+\frac{2\theta}{1+\theta}}_{q,1})}.
\end{align*}
\end{remark}
\medskip

We now turn to the proof of Theorem \ref{thm1.1}.
\medskip

\textbf{Proof of Theorem \ref{thm1.1}.}
Let us first recall the Leray projection operator to the divergence
free vector field space
$
\mathbb{P} \triangleq I-\nabla\Delta^{-1} \operatorname{div}.
$
Denote $\delta \triangleq d-\overline{d}_{0}$ and $\delta_{0}
\triangleq d_{0}-\overline{d}_{0}$. Then we can rewrite system
\eqref{eq1.1} as
\begin{align}\label{eq3.15}
\begin{cases}
&\partial_{t}u-\Delta u =-\mathbb{P}[u\cdot\nabla u
+\operatorname{div} (\nabla \delta\odot\nabla \delta)],\\
&\partial_{t}\delta-\Delta \delta=-u\cdot\nabla \delta+
|\nabla\delta|^{2}\delta
+|\nabla\delta|^{2}\overline{d}_{0},\\
&u|_{t=0}= u_{0}(x),\quad
\delta|_{t=0}=\delta_{0}(x),
\end{cases}
\end{align}
where the initial data satisfying
the following far field behavior
\begin{align}\label{eq3.16}
u_{0}\rightarrow 0, \delta_{0} \rightarrow 0 \text{ as } |x|\rightarrow
\infty.
\end{align}
By the Duhamel principle, we can express a solution $(u,\delta)$ of
\eqref{eq3.15} and \eqref{eq3.16} in the integral form:
\begin{align}\label{eq3.17}
\begin{cases}
&u(t) =e^{ t \Delta} u_{0}-\int_{0}^{t}
e^{(t-s)\Delta}\mathbb{P}[u\cdot\nabla u
+\operatorname{div} (\nabla \delta\odot\nabla \delta)](s)\text{d}s,\\
& \delta(t)= e^{
t\Delta}\delta_{0}+\int_{0}^{t}e^{(t-s)\Delta}[-u\cdot\nabla
\delta+ |\nabla\delta|^{2}\delta
+|\nabla\delta|^{2}\overline{d}_{0}](s)\text{d}s.
\end{cases}
\end{align}
Here, we notice that $\mathbb{P}$ is a homogeneous multiplier of
degree zero. Let
\begin{align*}
u_{L}\triangleq e^{ t\Delta} u_{0}  \text{ and }
\delta_{L}\triangleq e^{ t\Delta} \delta_{0}.
\end{align*}
Similar as the proof of Lemma \ref{lem3.1},
it is easy to see that $u_{L}\in
\widetilde{L}^{\infty}(\mathbb{R}_{+};e^{\sqrt{t}\Lambda_{1}}\dot{B}^{\frac{3}{p}-1}_{p,1}(\mathbb{R}^{3}))
\cap\widetilde{L}^{1}(\mathbb{R}_{+};e^{\sqrt{t}\Lambda_{1}}\dot{B}^{\frac{3}{p}+1}_{p,1}(\mathbb{R}^{3}))$
and $\delta_{L}\in
\widetilde{L}^{\infty}(\mathbb{R}_{+};e^{\sqrt{t}\Lambda_{1}}\dot{B}^{\frac{3}{q}}_{q,1}(\mathbb{R}^{3}))\cap
\widetilde{L}^{1}(\mathbb{R}_{+};e^{\sqrt{t}\Lambda_{1}}\dot{B}^{\frac{3}{q}+2}_{q,1}(\mathbb{R}^{3}))$,
and there hold
\begin{align}\label{eq3.18}
&\|u_{L}\|_{\widetilde{L}^{\infty}_{T}(e^{\sqrt{t}\Lambda_{1}}\dot{B}^{\frac{3}{p}-1}_{p,1})}
+
\|u_{L}\|_{\widetilde{L}^{1}_{T}(e^{\sqrt{t}\Lambda_{1}}\dot{B}^{\frac{3}{p}+1}_{p,1})}\leq
C_{0}\|u_{0}\|_{\dot{B}^{\frac{3}{p}-1}_{p,1}}, \text{ for all } 0<T\leq
\infty;\\
          \label{eq3.19}
&\|\delta_{L}\|_{\widetilde{L}^{\infty}_{T}(e^{\sqrt{t}\Lambda_{1}}\dot{B}^{\frac{3}{q}}_{q,1})}
+
\|\delta_{L}\|_{\widetilde{L}^{1}_{T}(e^{\sqrt{t}\Lambda_{1}}\dot{B}^{\frac{3}{q}+2}_{q,1})}\leq
C_{0}\|\delta_{0}\|_{\dot{B}^{\frac{3}{q}}_{q,1}}, \text{ for all }
0<T\leq \infty.
\end{align}
Define
$\overline{u}\triangleq u-u_{L}$ and $\overline{\delta}\triangleq
\delta-\delta_{L}$. Then $(u,\delta)$ is a solution of \eqref{eq3.17}
if and only if $(\overline{u},\overline{\delta})$ is the solution of
the following system
\begin{align}\label{eq3.20}
\begin{cases}
&\overline{u}(t) = -\int_{0}^{t}
e^{(t-s)\Delta}\mathbb{P}[\overline{u}\cdot\nabla \overline{u}
+\operatorname{div} (\nabla \overline{\delta}\odot\nabla \overline{\delta})+R_{1}+R_{2}](s)\text{d}s,\\
&
\overline{\delta}(t)=\int_{0}^{t}e^{(t-s)\Delta}[-\overline{u}\cdot\nabla
\overline{\delta}+ |\nabla\overline{\delta}|^{2}\overline{\delta}
+|\nabla\overline{\delta}|^{2}\overline{d}_{0}-R_{3}+ R_{4}+ R_{5}](s)\text{d}s,\\
\end{cases}
\end{align}
with
\begin{align*}
R_{1}& =u_{L}\cdot\nabla u_{L}+ u_{L}\cdot\nabla \overline{u}
+\overline{u}\cdot\nabla u_{L};\quad
R_{2} = \operatorname{div}(\nabla
\delta_{L}\odot\nabla\delta_{L})+2\operatorname{div}
(\nabla\delta_{L}\odot\nabla\overline{\delta});\nonumber\\
R_{3}&=u_{L}\cdot\nabla\delta_{L}+u_{L}\cdot\nabla\overline{\delta}
+\overline{u}\cdot\nabla\delta_{L};\qquad
R_{4}=|\nabla\delta_{L}|^{2}\overline{d}_{0}
+2|\nabla\delta_{L}||\nabla\overline{\delta}|\overline{d}_{0};\nonumber\\
R_{5}&=
|\nabla\delta_{L}|^{2}\delta_{L}+|\nabla\delta_{L}|^{2}\overline{\delta}
+2|\nabla\delta_{L}||\nabla\overline{\delta}|\delta_{L}+ 2
|\nabla\delta_{L}||\nabla\overline{\delta}|\overline{\delta}
+|\nabla\overline{\delta}|^{2}\delta_{L}.
\end{align*}

In what follows, we will use the iterative method to prove Theorem
\ref{thm1.1}. We only prove the local-in-time existence case, and in
a similar way, we can extend the local-in-time existence of solution
to the global-in-time existence of solution in the case of small
initial data.  Let $(\overline{u}_{0},\overline{\delta}_{0}):=
(0,0)$, and $(\overline{u}_{n}, \overline{\delta}_{n})$ satisfies
the following equations
\begin{align*}
\begin{cases}
&\!\!\!\overline{u}_{n}(t) = -\!\int_{0}^{t}
e^{(t-s)\Delta}\mathbb{P}[\overline{u}_{\!n-1}\cdot\!\nabla
\overline{u}_{\!n-1}\!
+\operatorname{div} (\nabla \overline{\delta}_{\!n-1}\odot\nabla \overline{\delta}_{\!n-1})\!+\!R_{1(n-1)}\!+\!R_{2(n-1)}](s)\text{d}s,\\
&
\!\!\!\overline{\delta}_{n}\!(t)\!=\!\!\int_{0}^{t}e^{(t-s)\Delta}[-\overline{u}_{n-1}\cdot\!\nabla
\overline{\delta}_{n-1}\!+\!
|\nabla\overline{\delta}_{n-1}|^{2}\overline{\delta}_{n-1}\!
+\!|\nabla\overline{\delta}_{n-1}|^{2}\overline{d}_{0}\!-\!R_{3(n-1)}\!+\!R_{4(n-1)}\!+\! R_{5(n-1)}](s)\text{d}s,\\
\end{cases}
\end{align*}
with
\begin{align*}
R_{1(n-1)}& =u_{L}\cdot\nabla u_{L}+ u_{L}\cdot\nabla
\overline{u}_{n-1}
+\overline{u}_{n-1}\cdot\nabla u_{L};\quad
R_{2(n-1)} = \operatorname{div}(\nabla
\delta_{L}\odot\nabla\delta_{L})+2\operatorname{div}
(\nabla\delta_{L}\odot\nabla\overline{\delta}_{n-1});\nonumber\\
R_{3(n-1)}&=u_{L}\cdot\nabla\delta_{L}+u_{L}\cdot\nabla\overline{\delta}_{n-1}
+\overline{u}_{n-1}\cdot\nabla\delta_{L}; \qquad
R_{4(n-1)}=|\nabla\delta_{L}|^{2}\overline{d}_{0}
+2|\nabla\delta_{L}||\nabla\overline{\delta}_{n-1}|\overline{d}_{0};\nonumber\\
R_{5(n-1)}&=
|\nabla\delta_{L}|^{2}\delta_{L}+|\nabla\delta_{L}|^{2}\overline{\delta}_{n-1}
+2|\nabla\delta_{L}||\nabla\overline{\delta}_{n-1}|\delta_{L}+ 2
|\nabla\delta_{L}||\nabla\overline{\delta}_{n-1}|\overline{\delta}_{n-1}
+|\nabla\overline{\delta}_{n-1}|^{2}\delta_{L}.
\end{align*}

\textbf{\textit{Step 1:} Uniform boundedness of $(\overline{u}_{n},\overline{\delta}_{n})$}

We claim that there exist $T>0$ and  small $\varepsilon>0$ such that
\begin{equation}\label{eq3.21}
  A\left(\|\overline{u}_{n}\|_{\widetilde{L}^{\infty}_{T}(e^{\sqrt{t}\Lambda_{1}}\dot{B}^{\frac{3}{p}-1}_{p,1})}
+
\|\overline{\delta}_{n}\|_{\widetilde{L}^{\infty}_{T}(e^{\sqrt{t}\Lambda_{1}}\dot{B}^{\frac{3}{q}}_{q,1})}\right)
+B\left(\|\overline{u}_{n}\|_{\widetilde{L}^{1}_{T}(e^{\sqrt{t}\Lambda_{1}}\dot{B}^{\frac{3}{p}+1}_{p,1})}
+\|\overline{\delta}_{n}\|_{\widetilde{L}^{1}_{T}(e^{\sqrt{t}\Lambda_{1}}\dot{B}^{\frac{3}{q}+2}_{q,1})}\right)\leq
\varepsilon,
\end{equation}
where
\begin{align}\label{eq3.22}
A\!\!=\!\!\left\{\begin{array}{ll}
\!\!\!  1, &\text{ if }
\!-\!\inf\{\frac{1}{3},\frac{1}{2p}\}\!\leq\!\frac{1}{q}\!-\!\frac{1}{p}\!<\!\frac{1}{3},\\
\!\!\! 2C_0C_{1}M_0,&\text{ if }
 \frac{1}{q}\!-\!\frac{1}{p}\!=\!\frac{1}{3},
\end{array}
\right.
\text{ and }
B\!\!=\!\!\left\{\begin{array}{ll}
 \!\!\! 1, &\text{ if }
\!-\!\inf\{\frac{1}{3},\frac{1}{2p}\}\!<\!\frac{1}{q}\!-\!\frac{1}{p}\!\leq\!\frac{1}{3},\\
\!\!\! 2C_0C_{1}M_0,&\text{ if }
 \frac{1}{q}\!-\!\frac{1}{p}\!=\!-\!\inf\{\frac{1}{3},\frac{1}{2p}\},
\end{array}
\right.
\end{align}
for all $n\in \mathbb{N}$. Here,
$M_{0}:=\|u_{0}\|_{\dot{B}^{\frac{3}{p}-1}_{p,1}}+\|\delta_{0}\|_{\dot{B}^{\frac{3}{p}}_{p,1}}$
is defined in Theorem \ref{thm1.1}, $C_{0}$ and $C_{1}$ are two uniform constants.
In fact, notice that $(\overline{u}_{0},\overline{\delta}_{0})=(0,0)$, it is easy to see \eqref{eq3.21}
holds for $n=0$.

 We are in a position to prove \eqref{eq3.21} for all $n\in\mathbb{N}_{+}$. In what follows,
  we shall only give the proof when $\frac{1}{q}-\frac{1}{p}=-\inf\{\frac{1}{3},\frac{1}{2p}\}$
 for the general $n\in\mathbb{N}_{+}$.  Due to
the proof  when $-\inf\{\frac{1}{3},\frac{1}{2p}\}<\frac{1}{q}-\frac{1}{p}\leq\frac{1}{3}$ is  similar,
 thus we omit the details.
By using Lemmas \ref{lem3.3}, \ref{lem3.5}, Remarks \ref{rem3.4}, \ref{rem3.6}, the
fact that the operator $\mathbb{P}$ is a homogeneous multiplier of
degree zero, and the interpolation inequality in Besov space,
we have  for some $\theta \in (0,1)$,
\begin{align}\label{eq3.23}
&\|\overline{u}_{n}\|_{\widetilde{L}^{\infty}_{T}(e^{\sqrt{t}\Lambda_{1}}\dot{B}^{\frac{3}{p}-1}_{p,1})}
+
\|\overline{u}_{n}\|_{\widetilde{L}^{1}_{T}(e^{\sqrt{t}\Lambda_{1}}\dot{B}^{\frac{3}{p}+1}_{p,1})}\nonumber\\
\leq&\!C\!\Big(\|\overline{u}_{n\!-\!1}\!\cdot\!\nabla \overline{u}_{n\!-\!1}
\|_{\widetilde{L}^{1}_{T}(e^{\sqrt{t}\Lambda_{1}}\dot{B}^{\frac{3}{p}-1}_{p,1})}
\!\!+\!\|\operatorname{div}(\nabla\overline{\delta}_{n\!-\!1}\!\odot\!\nabla
\overline{\delta}_{n\!-\!1})
\|_{\widetilde{L}^{1}_{T}(e^{\sqrt{t}\Lambda_{1}}\dot{B}^{\frac{3}{p}-1}_{p,1})}
\!\!+\!
\|R_{1(n\!-\!1)} \!+\!R_{2(n\!-\!1)}
\|_{\widetilde{L}^{1}_{T}(e^{\sqrt{t}\Lambda_{1}}\dot{B}^{\frac{3}{p}-1}_{p,1})}\!\Big)\nonumber\\
\leq &C\Big[
\Big(\|\overline{u}_{n\!-\!1}\|_{\widetilde{L}^{\infty}_{T}(e^{\sqrt{t}\Lambda_{1}}\dot{B}^{\frac{3}{p}-1}_{p,1})\cap
\widetilde{L}^{1}_{T}(e^{\sqrt{t}\Lambda_{1}}\dot{B}^{\frac{3}{p}+1}_{p,1})} +
\|\overline{\delta}_{n\!-\!1}\|_{\widetilde{L}^{\infty}_{T}(e^{\sqrt{t}\Lambda_{1}}\dot{B}^{\frac{3}{q}}_{q,1})\cap
\widetilde{L}^{1}_{T}(e^{\sqrt{t}\Lambda_{1}}\dot{B}^{\frac{3}{q}+2}_{q,1})}\Big)^{2}
\nonumber\\
&+\Big(\|u_{L}\|_{\widetilde{L}^{\infty}_{T}(e^{\sqrt{t}\Lambda_{1}}\dot{B}^{\frac{3}{p}-1}_{p,1})}
+\|\delta_{L}\|_{\widetilde{L}^{\infty}_{T}(e^{\sqrt{t}\Lambda_{1}}\dot{B}^{\frac{3}{q}}_{q,1})}\Big)
\Big(\|u_{L}\|_{\widetilde{L}^{1}_{T}(e^{\sqrt{t}\Lambda_{1}}\dot{B}^{\frac{3}{p}+1}_{p,1})}+
\|\delta_{L}\|_{\widetilde{L}^{1}_{T}(e^{\sqrt{t}\Lambda_{1}}\dot{B}^{\frac{3}{q}+2}_{q,1})}\Big)\nonumber\\
&+
\|u_{L}\|_{\widetilde{L}^{2}_{T}(e^{\sqrt{t}\Lambda_{1}}\dot{B}^{\frac{3}{p}}_{p,1})}
\Big(\|\overline{u}_{n\!-\!1}\|_{\widetilde{L}^{\infty}_{T}(e^{\sqrt{t}\Lambda_{1}}\dot{B}^{\frac{3}{p}-1}_{p,1})
\cap \widetilde{L}^{1}_{T}(e^{\sqrt{t}\Lambda_{1}}\dot{B}^{\frac{3}{p}+1}_{p,1})}
+\|\overline{\delta}_{n\!-\!1}\|_{\widetilde{L}^{\infty}_{T}(e^{\sqrt{t}\Lambda_{1}}\dot{B}^{\frac{3}{q}}_{q,1})
\cap \widetilde{L}^{1}_{T}(e^{\sqrt{t}\Lambda_{1}}\dot{B}^{\frac{3}{q}+2}_{q,1}) } \Big)\nonumber\\
&
+\|\delta_{L}\|_{\widetilde{L}^{1}_{T}(e^{\sqrt{t}\Lambda_{1}}\dot{B}^{\frac{3}{q}+2}_{q,1})}
\|\overline{\delta}_{n\!-\!1}\|_{\widetilde{L}^{\infty}_{T}(e^{\sqrt{t}\Lambda_{1}}\dot{B}^{\frac{3}{q}}_{q,1})}
 \Big]
+C_{1}
\|\delta_{L}\|_{\widetilde{L}^{\infty}_{T}(e^{\sqrt{t}\Lambda_{1}}\dot{B}^{\frac{3}{q}}_{q,1})}
\|\overline{\delta}_{n\!-\!1}\|_{\widetilde{L}^{1}_{T}(e^{\sqrt{t}\Lambda_{1}}\dot{B}^{\frac{3}{q}+2}_{q,1})},
\end{align}
and
\begin{align}\label{eq3.24}
&\|\overline{\delta}_{n}\|_{\widetilde{L}^{\infty}_{T}(e^{\sqrt{t}\Lambda_{1}}\dot{B}^{\frac{3}{q}}_{q,1})}
+
\|\overline{\delta}_{n}\|_{\widetilde{L}^{1}_{T}(e^{\sqrt{t}\Lambda_{1}}\dot{B}^{\frac{3}{q}+2}_{q,1})}\nonumber\\
\leq& C\Big(\|\overline{u}_{n\!-\!1}\cdot\nabla \overline{\delta}_{n\!-\!1}
\|_{\widetilde{L}^{1}_{T}(e^{\sqrt{t}\Lambda_{1}}\dot{B}^{\frac{3}{q}}_{q,1})}
+\||\nabla\overline{\delta}_{n\!-\!1}|^{2}\overline{\delta}_{n\!-\!1}
\|_{\widetilde{L}^{1}_{T}(e^{\sqrt{t}\Lambda_{1}}\dot{B}^{\frac{3}{q}}_{q,1})}
+\||\nabla\overline{\delta}_{n\!-\!1}|^{2}\overline{d}_{0}
\|_{\widetilde{L}^{1}_{T}(e^{\sqrt{t}\Lambda_{1}}\dot{B}^{\frac{3}{q}}_{q,1})}\nonumber\\
 &+
\|-R_{3(n\!-\!1)}+ R_{4(n\!-\!1)}+ R_{5(n\!-\!1)}
\|_{\widetilde{L}^{1}_{T}(e^{\sqrt{t}\Lambda_{1}}\dot{B}^{\frac{3}{q}}_{q,1})}
\Big)\nonumber\\
\leq&\!C\Big[\!
(1\!+\!\!\|\delta_{L}\|_{\widetilde{L}^{\infty}_{T}(e^{\sqrt{t}\Lambda_{1}}\dot{B}^{\frac{3}{q}}_{q,1})})
\Big(\|\overline{u}_{n\!-\!1}\|_{\widetilde{L}^{\infty}_{T}(e^{\sqrt{t}\Lambda_{1}}\dot{B}^{\frac{3}{p}-1}_{p,1})\cap
\widetilde{L}^{1}_{T}(e^{\sqrt{t}\Lambda_{1}}\dot{B}^{\frac{3}{p}+1}_{p,1})}\!\!+\!
\|\overline{\delta}_{n\!-\!1}\|_{\widetilde{L}^{\infty}_{T}(e^{\sqrt{t}\Lambda_{1}}\dot{B}^{\frac{3}{q}}_{q,1})\cap
\widetilde{L}^{1}_{T}(e^{\sqrt{t}\Lambda_{1}}\dot{B}^{\frac{3}{q}+2}_{q,1})}\Big)^{2}
\nonumber\\
&+\Big(\|u_{L}\|_{\widetilde{L}^{\infty}_{T}(e^{\sqrt{t}\Lambda_{1}}\dot{B}^{\frac{3}{p}-1}_{p,1})}
+\|\delta_{L}\|_{\widetilde{L}^{\infty}_{T}(e^{\sqrt{t}\Lambda_{1}}\dot{B}^{\frac{3}{q}}_{q,1})}\Big)
\Big(\|u_{L}\|_{\widetilde{L}^{1}_{T}(e^{\sqrt{t}\Lambda_{1}}\dot{B}^{\frac{3}{p}+1}_{p,1})}+
\|\delta_{L}\|_{\widetilde{L}^{1}_{T}(e^{\sqrt{t}\Lambda_{1}}\dot{B}^{\frac{3}{q}+2}_{q,1})}\Big)
\nonumber\\
&
+\|\overline{\delta}_{n\!-\!1}\|^{2}_{\widetilde{L}^{\infty}_{T}(e^{\sqrt{t}\Lambda_{1}}\dot{B}^{\frac{3}{q}}_{q,1})}
\|\overline{\delta}_{n\!-\!1}\|_{\widetilde{L}^{1}_{T}(e^{\sqrt{t}\Lambda_{1}}\dot{B}^{\frac{3}{q}+2}_{q,1})}
+\|\delta_{L}\|^{2}_{\widetilde{L}^{\infty}_{T}(e^{\sqrt{t}\Lambda_{1}}\dot{B}^{\frac{3}{q}}_{q,1})}
\|\delta_{L}\|_{\widetilde{L}^{1}_{T}(e^{\sqrt{t}\Lambda_{1}}\dot{B}^{\frac{3}{q}+2}_{q,1})}\nonumber\\
&+\Big(\|u_{L}\|_{\widetilde{L}^{2}_{T}(e^{\sqrt{t}\Lambda_{1}}\dot{B}^{\frac{3}{p}}_{p,1})}
+\|\delta_{L}\|_{\widetilde{L}^{2}_{T}(e^{\sqrt{t}\Lambda_{1}}\dot{B}^{\frac{3}{q}+1}_{q,1})}
+\|u_{L}\|_{\widetilde{L}^{1+\theta}_{T}
(e^{\sqrt{t}\Lambda_{1}}\dot{B}^{\frac{3}{p}+\frac{1-\theta}{1+\theta}}_{p,1})}
+\|\delta_{L}\|_{\widetilde{L}^{\frac{1+\theta}{\theta}}_{T}
(e^{\sqrt{t}\Lambda_{1}}\dot{B}^{\frac{3}{q}+\frac{2\theta}{1+\theta}}_{q,1})}\Big)
\nonumber\\
&\times \Big(\|\overline{u}_{n\!-\!1}\|_{\widetilde{L}^{\infty}_{T}(e^{\sqrt{t}\Lambda_{1}}\dot{B}^{\frac{3}{p}-1}_{p,1})
\cap \widetilde{L}^{1}_{T}(e^{\sqrt{t}\Lambda_{1}}\dot{B}^{\frac{3}{p}+1}_{p,1})}
+\|\overline{\delta}_{n\!-\!1}\|_{\widetilde{L}^{\infty}_{T}(e^{\sqrt{t}\Lambda_{1}}\dot{B}^{\frac{3}{q}}_{q,1})
\cap \widetilde{L}^{1}_{T}(e^{\sqrt{t}\Lambda_{1}}\dot{B}^{\frac{3}{q}+2}_{q,1}) } \Big)\Big],
\end{align}
where we have used the fact that
$\|\nabla f\|_{L^{r}_{T}(e^{\sqrt{t}\Lambda_{1}}\dot{B}^{s}_{p,q})}\approx
\|f\|_{L^{r}_{T}(e^{\sqrt{t}\Lambda_{1}}\dot{B}^{s+1}_{p,q})}$ with $s\in\mathbb{R},1\leq p,q,r\leq \infty$.
On the
other hand, notice that from Lemma \ref{lem3.1}, we find that there
exists $T_{0}>0$ such that for any given small number $0<\zeta<1$,
it holds that for all $\theta\in (0,1]$ and $0<T\leq T_{0}$,
\begin{align} \label{eq3.25}
&\|u_{L}\|_{\widetilde{L}^{1}_{T}(e^{\sqrt{t}\Lambda_{1}}\dot{B}^{\frac{3}{p}+1}_{p,1})}
+\|\delta_{L}\|_{\widetilde{L}^{1}_{T}(e^{\sqrt{t}\Lambda_{1}}\dot{B}^{\frac{3}{q}+2}_{q,1})} \leq
\zeta,  \\
     \label{eq3.26}
&\|u_{L}\|_{\widetilde{L}^{1+\theta}_{T}
(e^{\sqrt{t}\Lambda_{1}}\dot{B}^{\frac{3}{p}+\frac{1-\theta}{1+\theta}}_{p,1})\cap
\widetilde{L}^{\frac{1+\theta}{\theta}}_{T}
(e^{\sqrt{t}\Lambda_{1}}\dot{B}^{\frac{3}{p}+\frac{\theta-1}{1+\theta}}_{p,1})}
+\|\delta_{L}\|_{\widetilde{L}^{1+\theta}_{T}
(e^{\sqrt{t}\Lambda_{1}}\dot{B}^{\frac{3}{q}+\frac{2}{1+\theta}}_{q,1})\cap
\widetilde{L}^{\frac{1+\theta}{\theta}}_{T}
(e^{\sqrt{t}\Lambda_{1}}\dot{B}^{\frac{3}{q}+\frac{2\theta}{1+\theta}}_{q,1})}\leq
\zeta.
\end{align}
From inequalities \eqref{eq3.18} and \eqref{eq3.19}, we have
\begin{align}\label{eq3.27}
\|u_{L}\|_{\widetilde{L}^{\infty}_{T}(e^{\sqrt{t}\Lambda_{1}}\dot{B}^{\frac{3}{p}-1}_{p,1})}
+\|\delta_{L}\|_{\widetilde{L}^{\infty}_{T}(e^{\sqrt{t}\Lambda_{1}}\dot{B}^{\frac{3}{q}}_{q,1})}
\leq
C_{0}(\|u_{0}\|_{\dot{B}^{\frac{3}{p}-1}_{p,1}}+\|\delta_{0}\|_{\dot{B}^{\frac{3}{q}}_{q,1}})
=C_{0}M_{0}.
\end{align}
Hence, by using the fact that $(\overline{u}_{n\!-\!1}(0),
\overline{\delta}_{n\!-\!1}(0))= (0,0)$, we can choose $T\leq T_{0}$
small enough such that the claim \eqref{eq3.21} holds for $n\!-\!1$,
i.e.,
\begin{eqnarray}\label{eq3.28}
\|\overline{u}_{n\!-\!1}\|_{\widetilde{L}^{\infty}_{T}(e^{\sqrt{t}\Lambda_{1}}\dot{B}^{\frac{3}{p}-1}_{p,1})
\cap\widetilde{L}^{1}_{T}(e^{\sqrt{t}\Lambda_{1}}\dot{B}^{\frac{3}{p}+1}_{p,1})}
+2C_{0}C_{1}M_{0}
\|\overline{\delta}_{n\!-\!1}\|_{\widetilde{L}^{\infty}_{T}(e^{\sqrt{t}\Lambda_{1}}\dot{B}^{\frac{3}{q}}_{q,1})
\cap\widetilde{L}^{1}_{T}(e^{\sqrt{t}\Lambda_{1}}\dot{B}^{\frac{3}{q}+2}_{q,1})}\leq
\varepsilon.
\end{eqnarray}
Inserting \eqref{eq3.25}--\eqref{eq3.28} into \eqref{eq3.23} and
\eqref{eq3.24}, it follows that
\begin{align}\label{eq3.29}
&\|\overline{u}_{n}\|_{\widetilde{L}^{\infty}_{T}(e^{\sqrt{t}\Lambda_{1}}\dot{B}^{\frac{3}{p}-1}_{p,1})
\cap \widetilde{L}^{1}_{T}(e^{\sqrt{t}\Lambda_{1}}\dot{B}^{\frac{3}{p}+1}_{p,1})}
+2C_{0}C_{1}M_{0}
\|\overline{\delta}_{n}\|_{\widetilde{L}^{\infty}_{T}(e^{\sqrt{t}\Lambda_{1}}\dot{B}^{\frac{3}{q}}_{q,1})
\cap \widetilde{L}^{1}_{T}(e^{\sqrt{t}\Lambda_{1}}\dot{B}^{\frac{3}{q}+2}_{q,1})}\nonumber\\
\leq &C\Big[
(1+\|\delta_{L}\|_{\widetilde{L}^{\infty}_{T}(e^{\sqrt{t}\Lambda_{1}}\dot{B}^{\frac{3}{q}}_{q,1})})
\Big(\|\overline{u}_{n\!-\!1}\|_{\widetilde{L}^{\infty}_{T}(e^{\sqrt{t}\Lambda_{1}}\dot{B}^{\frac{3}{p}-1}_{p,1})\cap
\widetilde{L}^{1}_{T}(e^{\sqrt{t}\Lambda_{1}}\dot{B}^{\frac{3}{p}+1}_{p,1})}+
 \|\overline{\delta}_{n\!-\!1}\|_{\widetilde{L}^{\infty}_{T}(e^{\sqrt{t}\Lambda_{1}}\dot{B}^{\frac{3}{q}}_{q,1})
 \cap
\widetilde{L}^{1}_{T}(e^{\sqrt{t}\Lambda_{1}}\dot{B}^{\frac{3}{q}+2}_{q,1})}\Big)^{2}
\nonumber\\
&+\Big(\|u_{L}\|_{\widetilde{L}^{\infty}_{T}(e^{\sqrt{t}\Lambda_{1}}\dot{B}^{\frac{3}{p}-1}_{p,1})}
+\|\delta_{L}\|_{\widetilde{L}^{\infty}_{T}(e^{\sqrt{t}\Lambda_{1}}\dot{B}^{\frac{3}{q}}_{q,1})}\Big)
\Big(\|u_{L}\|_{\widetilde{L}^{1}_{T}(e^{\sqrt{t}\Lambda_{1}}\dot{B}^{\frac{3}{p}+1}_{p,1})}+
\|u_{L}\|_{\widetilde{L}^{1}_{T}(e^{\sqrt{t}\Lambda_{1}}\dot{B}^{\frac{3}{p}+2}_{p,1})}\Big)
\nonumber\\
&+\Big(\|u_{L}\|_{\widetilde{L}^{2}_{T}(e^{\sqrt{t}\Lambda_{1}}\dot{B}^{\frac{3}{p}}_{p,1})}+\|u_{L}\|_{\widetilde{L}^{1+\theta}_{T}(e^{\sqrt{t}\Lambda_{1}}\dot{B}^{\frac{3}{p}+\frac{1-\theta}{1+\theta}}_{p,1})}
+\|\delta_{L}\|_{\widetilde{L}^{2}_{T}(e^{\sqrt{t}\Lambda_{1}}\dot{B}^{\frac{3}{q}+1}_{q,1})}
+\|\delta_{L}\|_{\widetilde{L}^{\frac{1+\theta}{\theta}}_{T}(e^{\sqrt{t}\Lambda_{1}}
\dot{B}^{\frac{3}{q}+\frac{2\theta}{1+\theta}}_{q,1})}\Big)
\nonumber\\
&\times
\Big(\|\overline{u}_{n\!-\!1}\|_{\widetilde{L}^{\infty}_{T}(e^{\sqrt{t}\Lambda_{1}}\dot{B}^{\frac{3}{p}-1}_{p,1})\cap
\widetilde{L}^{1}_{T}(e^{\sqrt{t}\Lambda_{1}}\dot{B}^{\frac{3}{p}+1}_{p,1})}+
 \|\overline{\delta}_{n\!-\!1}\|_{ \widetilde{L}^{\infty}_{T}(e^{\sqrt{t}\Lambda_{1}}\dot{B}^{\frac{3}{q}}_{ q,1})
 \cap
\widetilde{L}^{1}_{T}(e^{\sqrt{t}\Lambda_{1}}\dot{B}^{ \frac{3}{q}+2}_{q,1})} \Big) \nonumber\\
&
+ \|\overline{\delta}_{n\!-\!1}\|^{2}_{ \widetilde{L}^{\infty}_{T}
(e^{\sqrt{t}\Lambda_{1}}\dot{B}^{\frac{3}{q}}_{ q,1})}
 \|\overline{\delta}_{n\!-\!1}\|_{ \widetilde{L}^{1}_{T}(e^{\sqrt{t}\Lambda_{1}}\dot{B}^{ \frac{3}{q}+2}_{q,1})}
+ \|\delta_{L}\|^{2}_{
\widetilde{L}^{\infty}_{T}(e^{\sqrt{t}\Lambda_{1}}\dot{B}^{\frac{3}{q}}_{ q,1})}
 \|\delta_{L}\|_{ \widetilde{L}^{1}_{T}(e^{\sqrt{t}\Lambda_{1}}\dot{B}^{ \frac{3}{q}+2}_{q,1})} \nonumber\\
&
 +\|\delta_{L}\|_{\widetilde{L}^{1}_{T}(e^{\sqrt{t}\Lambda_{1}}\dot{B}^{\frac{3}{q}+2}_{q,1})}
\|\overline{\delta}_{n\!-\!1}\|_{\widetilde{L}^{\infty}_{T}(e^{\sqrt{t}\Lambda_{1}}\dot{B}^{\frac{3}{q}}_{q,1})}
\Big]
+C_{0}C_{1}M_{0}
\|\overline{\delta}_{n\!-\!1}\|_{\widetilde{L}^{\infty}_{T}(e^{\sqrt{t}\Lambda_{1}}\dot{B}^{\frac{3}{q}+2}_{q,1})}\nonumber\\
\leq& C\Big[(1+C_{0}M_{0})
(\varepsilon+\frac{\varepsilon}{2C_{0}C_{1}M_{0}})^{2}
+\zeta(C_{0}M_{0}+(C_{0}M_{0})^{2}+2\varepsilon+\frac{\varepsilon}{2C_{0}C_{1}M_{0}})
+\frac{\varepsilon^{3}}{(2C_{0}C_{1}M_{0})^{3}}\Big]+\frac{\varepsilon}{2}.
\end{align}
By letting $\zeta,\varepsilon>0$ small enough such that,
\begin{align*}
&C(1+C_{0}M_{0})
(1+\frac{1}{2C_{0}C_{1}M_{0}})^{2}\varepsilon\leq
\frac{1}{8};\qquad
C(C_{0}M_{0}+C_{0}^{2}M_{0}^{2})\zeta\leq \frac{\varepsilon}{8};\nonumber\\
&C\zeta(2+\frac{1}{2C_{0}C_{1}M_{0}})
+\frac{\varepsilon^{2}\mu^{2}}{(2C_{0}C_{1}M_{0})^{3}}\leq
\frac{1}{4}.
\end{align*}
Then, it follows from \eqref{eq3.29} that
\begin{align*}
\|\overline{u}_{n}\|_{\widetilde{L}^{\infty}_{T}(e^{\sqrt{t}\Lambda_{1}}\dot{B}^{\frac{3}{p}-1}_{p,1})
\cap\widetilde{L}^{1}_{T}(e^{\sqrt{t}\Lambda_{1}}\dot{B}^{\frac{3}{p}+1}_{p,1})}
+2C_{0}C_{1}M_{0}
\|\overline{\delta}_{n}\|_{\widetilde{L}^{\infty}_{T}(e^{\sqrt{t}\Lambda_{1}}\dot{B}^{\frac{3}{q}}_{q,1})
\cap\widetilde{L}^{1}_{T}(e^{\sqrt{t}\Lambda_{1}}\dot{B}^{\frac{3}{q}+2}_{q,1})}\leq
\varepsilon.
\end{align*}
 Arguing by induction, we conclude that \eqref{eq3.21} with $\frac{1}{q}-\frac{1}{p}=-\inf\{\frac{1}{3},\frac{1}{2p}\}$
holds for all $n\in \mathbb{N}_{+}$.
\medskip

\textbf{\textit{Step 2:} Convergence of $(\overline{u}_{n},\overline{\delta}_{n})$}

We claim that $(\overline{u}_{n},\overline{\delta}_{n})$ is a Cauchy
sequence in
$\widetilde{L}^{\infty}_{T}(e^{\sqrt{t}\Lambda_{1}}\dot{B}^{\frac{3}{p}-1}_{p,1}\!(\mathbb{R}^{3}))\cap
\widetilde{L}^{1}_{T}(e^{\sqrt{t}\Lambda_{1}}\dot{B}^{\frac{3}{p}+1}_{p,1}\!(\mathbb{R}^{3}))\times
\widetilde{L}^{\infty}_{T}(e^{\sqrt{t}\Lambda_{1}}\dot{B}^{\frac{3}{q}}_{q,1}\!(\mathbb{R}^{3}))\cap
\widetilde{L}^{1}_{T}(e^{\sqrt{t}\Lambda_{1}}\dot{B}^{\frac{3}{q}+2}_{q,1}\!(\mathbb{R}^{3}))$. Let us
consider
\begin{align}\label{eq3.30}
&A\Big(\|\overline{u}_{m+n+1}-\overline{u}_{n+1}\|_{
\widetilde{L}^{\infty}_{T}(e^{\sqrt{t}\Lambda_{1}}\dot{B}^{\frac{3}{p}-1}_{p,1})
\cap \widetilde{L}^{1}_{T}(e^{\sqrt{t}\Lambda_{1}}\dot{B}^{\frac{3}{p}+1}_{p,1})}\Big)
+B
 \|\overline{\delta}_{m+n+1}-\overline{\delta}_{n+1}\|_{
 \widetilde{L}^{\infty}_{T}(e^{\sqrt{t}\Lambda_{1}}\dot{B}^{\frac{3}{q}}_{q,1})
 \cap\widetilde{L}^{1}_{T}(e^{\sqrt{t}\Lambda_{1}}\dot{B}^{\frac{3}{q}+2}_{q,1})},
\end{align}
for all $m,n\in\mathbb{N}_{+}$, where $A, B$ defined as \eqref{eq3.22}. Similar as the previous subsection,
we only consider the case of  $\frac{1}{q}-\frac{1}{p}=-\inf\{\frac{1}{3},\frac{1}{2p}\}$, and in the similar
way, we can establish that
$(\overline{u}_{n},\overline{\delta}_{n})$ is still a Cauchy
sequence in the case of $-\inf\{\frac{1}{3},\frac{1}{2p}\}<\frac{1}{q}-\frac{1}{p}\leq\frac{1}{3}$. According
to the proof of Lemma \ref{lem3.2}, we need
to estimate
$\widetilde{L}^{1}_{T}(e^{\sqrt{t}\Lambda_{1}}\dot{B}^{\frac{3}{p}-1}_{p,1}(\mathbb{R}^{3}))$-norm or
$\widetilde{L}^{1}_{T}(e^{\sqrt{t}\Lambda_{1}}\dot{B}^{\frac{3}{q}}_{q,1}(\mathbb{R}^{3}))$-norm of the
following terms:
\begin{align*}
I_{1}(t,x)&\triangleq \int_{0}^{t} e^{(t-s)\Delta}
\mathbb{P}[(\overline{u}_{m+n}\cdot\nabla
\overline{u}_{m+n})-(\overline{u}_{n}\cdot\nabla
\overline{u}_{n})](s)\text{d}s;\\
I_{2}(t,x)&\triangleq \int_{0}^{t}e^{(t-s)\Delta} \mathbb{P}
[\nabla\cdot(\nabla\overline{\delta}_{m+n}\odot\nabla\overline{\delta}_{m+n})
-\nabla\cdot(\nabla\overline{\delta}_{n}\odot\nabla\overline{\delta}_{n})](s)\text{d}s;\\
I_{3}(t,x)&\triangleq\int_{0}^{t}e^{(t-s)\Delta}
[\overline{u}_{m+n}\cdot\nabla
\overline{\delta}_{m+n}-\overline{u}_{n}\cdot\nabla\overline{\delta}_{n}](s)\text{d}s;\\
I_{4}(t,x)&\triangleq \int_{0}^{t}e^{(t-s)\Delta}
[|\nabla\overline{\delta}_{m+n}|^{2}\overline{\delta}_{m+n}-|\nabla\overline{\delta}_{n}|^{2}\overline{\delta}_{n}](s)\text{d}s;\\
I_{5}(t,x)&\triangleq \int_{0}^{t} e^{(t-s)\Delta}
[|\nabla\overline{\delta}_{m+n}|^{2}\overline{d}_{0}-|\nabla\overline{\delta}_{n}|^{2}\overline{d}_{0}](s)\text{d}s;\\
J_{j}(t,x)&\triangleq \int_{0}^{t} e^{(t-s)\Delta} [R_{j(m+n)}-R_{j(n)}](s)\text{d}s\quad \text{ with } j=1,2,\cdots, 5.
\end{align*}
Then, similar as derivation of \eqref{eq3.21},  if we choose
$\varepsilon,\zeta$ small enough, a  straight calculus  by using
Lemmas \ref{lem3.3}, \ref{lem3.5}, Remarks \ref{rem3.4}, \ref{rem3.6}, the
interpolation inequality, and the inequalities
\eqref{eq3.19} and \eqref{eq3.25}--\eqref{eq3.27} gives that
\begin{align*}
&\|\overline{u}_{m+n+1}-\overline{u}_{n+1}\|_{
\widetilde{L}^{\infty}_{T}(e^{\sqrt{t}\Lambda_{1}}\dot{B}^{\frac{3}{p}-1}_{p,1})
\cap \widetilde{L}^{1}_{T}(e^{\sqrt{t}\Lambda_{1}}\dot{B}^{\frac{3}{p}+1}_{p,1})}
+
2C_{0}C_{1}M_{0}
\|\overline{\delta}_{m+n+1}-\overline{\delta}_{n+1}\|_{
\widetilde{L}^{\infty}_{T}(e^{\sqrt{t}\Lambda_{1}}\dot{B}^{\frac{3}{q}}_{q,1})
\cap\widetilde{L}^{1}_{T}(e^{\sqrt{t}\Lambda_{1}}\dot{B}^{\frac{3}{q}+2}_{q,1})}\nonumber\\
\leq&
\frac{1}{2}\Big(\|\overline{u}_{m+n}-\overline{u}_{n}\|_{
\widetilde{L}^{\infty}_{T}(e^{\sqrt{t}\Lambda_{1}}\dot{B}^{\frac{3}{p}-1}_{p,1})
\cap\widetilde{L}^{1}_{T}(e^{\sqrt{t}\Lambda_{1}}\dot{B}^{\frac{3}{p}+1}_{p,1})}
+
2C_{0}C_{1}M_{0}
\|\overline{\delta}_{m+n}-\overline{\delta}_{n}\|_{
\widetilde{L}^{\infty}_{T}(e^{\sqrt{t}\Lambda_{1}}\dot{B}^{\frac{3}{q}}_{q,1})
\cap \widetilde{L}^{1}_{T}(e^{\sqrt{t}\Lambda_{1}}\dot{B}^{\frac{3}{q}+2}_{q,1})}\Big),
\end{align*}
when $
\frac{1}{q}-\frac{1}{p}=-\inf\{\frac{1}{3},\frac{1}{2p}\}$. Thus we obtain \eqref{eq3.30} for all
$m,n\in\mathbb{N}_{+}$.

From \eqref{eq3.30}, we see  that $(u_{n}, \delta_{n})$ is a Cauchy sequence in
$\widetilde{L}^{\infty}_{T}(e^{\sqrt{t}\Lambda_{1}}\dot{B}^{\frac{3}{p}-1}_{p,1}(\mathbb{R}^{3}))\cap
\widetilde{L}^{1}_{T}(e^{\sqrt{t}\Lambda_{1}}\dot{B}^{\frac{3}{p}+1}_{p,1}(\mathbb{R}^{3}))\!\times
\widetilde{L}^{\infty}_{T}(e^{\sqrt{t}\Lambda_{1}}\dot{B}^{\frac{3}{q}}_{q,1}(\mathbb{R}^{3}))\cap
\widetilde{L}^{1}_{T}(e^{\sqrt{t}\Lambda_{1}}\dot{B}^{\frac{3}{q}+2}_{q,1}(\mathbb{R}^{3}))$ with $1\leq
p,q<\infty$ and  $-\inf\{\frac{1}{3},\frac{1}{2p}\}\leq
\frac{1}{q}-\frac{1}{p}\leq \frac{1}{3}$. Let
\begin{align*}
 \overline{u}_{n}\rightarrow \overline{u} \text{ in
 }\widetilde{L}^{\infty}_{T}(e^{\sqrt{t}\Lambda_{1}}\dot{B}^{\frac{3}{p}-1}_{p,1}(\mathbb{R}^{3}))\cap
\widetilde{L}^{1}_{T}(e^{\sqrt{t}\Lambda_{1}}\dot{B}^{\frac{3}{p}+1}_{p,1}(\mathbb{R}^{3})),
\end{align*}
 and
 \begin{align*}
  \overline{\delta}_{n}\rightarrow \overline{\delta} \text{ in
 }\widetilde{L}^{\infty}_{T}(e^{\sqrt{t}\Lambda_{1}}\dot{B}^{\frac{3}{q}}_{q,1}(\mathbb{R}^{3}))
 \cap\widetilde{L}^{1}_{T}(e^{\sqrt{t}\Lambda_{1}}\dot{B}^{\frac{3}{q}+2}_{q,1}(\mathbb{R}^{3})).
\end{align*}
By using the standard arguments, we can prove that
$(\overline{u},\overline{\delta})$ is a solution of \eqref{eq3.20} on
$\mathbb{R}^{3}\times[0,T]$.  This implies that
\begin{align*}
(u,\delta)\!=\!
(u_{\!L}\!+\!\overline{u}, \delta_{\!L}\!+\!\overline{\delta})\!
\in \!\widetilde{L}^{\infty}_{T}(e^{\!\sqrt{t}\Lambda_{1}}\dot{B}^{\frac{3}{p}-1}_{p,1}\!(\mathbb{R}^{3}))
\!\cap\! \widetilde{L}^{1}_{T}(e^{\!\sqrt{t}\Lambda_{1}}\dot{B}^{\frac{3}{p}+1}_{p,1}\!(\mathbb{R}^{3}))
\!\times\! \widetilde{L}^{\infty}_{T}(e^{\!\sqrt{t}\Lambda_{1}}\dot{B}^{\frac{3}{q}}_{q,1}\!(\mathbb{R}^{3}))
\!\cap\! \widetilde{L}^{1}_{T}(e^{\!\sqrt{t}\Lambda_{1}}\dot{B}^{\frac{3}{q}+2}_{q,1}\!(\mathbb{R}^{3}))
\end{align*}
 is the solution of \eqref{eq3.17}. In fact, the solution can be extended step by step  and finally  we
 have a maximal time $T^{*}$ verifying $(u,\delta)\in \widetilde{L}^{\infty}_{T^{*}}
 (e^{\sqrt{t}\Lambda_{1}}\dot{B}^{\frac{3}{p}-1}_{p,1}(\mathbb{R}^{3}))
\cap \widetilde{L}^{1}_{T^{*}}(e^{\sqrt{t}\Lambda_{1}}\dot{B}^{\frac{3}{p}+1}_{p,1}(\mathbb{R}^{3}))
\times \widetilde{L}^{\infty}_{T^{*}}(e^{\sqrt{t}\Lambda_{1}}\dot{B}^{\frac{3}{q}}_{q,1}(\mathbb{R}^{3}))
\cap \widetilde{L}^{1}_{T^{*}}(e^{\sqrt{t}\Lambda_{1}}\dot{B}^{\frac{3}{q}+2}_{q,1}(\mathbb{R}^{3}))$.
If $T^{*}<\infty$ and for all
$\theta\in(0,1]$
\begin{align*}
&\|u\|_{\widetilde{L}^{1}_{T^{*}}( e^{\sqrt{t}\Lambda_{1}}\dot{B}^{\frac{3}{p}+1}_{p,1})}
+\|u\|_{\widetilde{L}^{1+\theta}_{T^{*}}
(e^{\sqrt{t}\Lambda_{1}}\dot{B}^{\frac{3}{p}+\frac{1-\theta}{1+\theta}}_{p,1})\cap
\widetilde{L}^{\frac{1+\theta}{\theta}}_{T^{*}}
(e^{\sqrt{t}\Lambda_{1}}\dot{B}^{\frac{3}{p}+\frac{\theta-1}{1+\theta}}_{p,1})}\nonumber\\
&+\|\delta\|_{\widetilde{L}^{1}_{T^{*}}(  e^{\sqrt{t}\Lambda_{1}}\dot{B}^{\frac{3}{q}+2}_{q,1})}
+ \|\delta\|_{\widetilde{L}^{1+\theta}_{T^{*}}
(e^{\sqrt{t}\Lambda_{1}}\dot{B}^{\frac{3}{q}+\frac{2}{1+\theta}}_{q,1})\cap
\widetilde{L}^{\frac{1+\theta}{\theta}}_{T^{*}}
(e^{\sqrt{t}\Lambda_{1}}\dot{B}^{\frac{3}{q}+\frac{2\theta}{1+\theta}}_{q,1})}
<\infty,
\end{align*}
 we claim that
the solution $(u,d)$ can be extended beyond $T^{*}$. Indeed, let us consider the integral equations
\begin{align*}
\begin{cases}
&u(t) =e^{ (t-T) \Delta} u(T)-\int_{T}^{t}
e^{(t-s)\Delta}\mathbb{P}[u\cdot\nabla u
+\operatorname{div} (\nabla \delta\odot\nabla \delta)](s)\text{d}s,\\
& \delta(t)= e^{
(t-T)\Delta}\delta(T)+\int_{T}^{t}e^{(t-s)\Delta}[-u\cdot\nabla
\delta+ |\nabla\delta|^{2}\delta
+|\nabla\delta|^{2}\overline{d}_{0}](s)\text{d}s.
\end{cases}
\end{align*}
Let $\mathcal{A}(u,\delta):=\mathbb{P}[u\cdot\nabla u
+\operatorname{div} (\nabla \delta\odot\nabla \delta)]$
 and $\mathcal{B}(u,\delta):=-u\cdot\nabla
\delta+ |\nabla\delta|^{2}\delta
+|\nabla\delta|^{2}\overline{d}_{0}$. By using Lemma \ref{lem3.2},  we see that
\begin{align}\label{eq3.31}
&\|e^{ (t-T) \Delta} u(T)\|_{\widetilde{L}^{1}(T,T^{*};e^{\sqrt{t}\Lambda_{1}}\dot{B}^{\frac{3}{p}+1}_{p,1})}
+\|e^{ (t-T) \Delta} u(T)\|_{\widetilde{L}^{1+\theta}(T,T^{*};
(e^{\sqrt{t}\Lambda_{1}}\dot{B}^{\frac{3}{p}+\frac{1-\theta}{1+\theta}}_{p,1})\cap
\widetilde{L}^{\frac{1+\theta}{\theta}}(T,T^{*};
(e^{\sqrt{t}\Lambda_{1}}\dot{B}^{\frac{3}{p}+\frac{\theta-1}{1+\theta}}_{p,1})}\nonumber\\
&+ \|e^{ (t-T) \Delta} \delta(T)\|_{\widetilde{L}^{1}(T,T^{*};
e^{\sqrt{t}\Lambda_{1}}\dot{B}^{\frac{3}{q}+2}_{q,1})}
+ \|e^{ (t-T) \Delta} \delta(T)\|_{\widetilde{L}^{1+\theta}(T,T^{*};
(e^{\sqrt{t}\Lambda_{1}}\dot{B}^{\frac{3}{q}+\frac{2}{1+\theta}}_{q,1})\cap
\widetilde{L}^{\frac{1+\theta}{\theta}}(T,T^{*};
(e^{\sqrt{t}\Lambda_{1}}\dot{B}^{\frac{3}{q}+\frac{2\theta}{1+\theta}}_{q,1})}\nonumber\\
\leq& \|u\|_{\widetilde{L}^{1}(T,T^{*};e^{\sqrt{t}\Lambda_{1}}\dot{B}^{\frac{3}{p}+1}_{p,1})}
+\|u\|_{\widetilde{L}^{1+\theta}(T,T^{*};
(e^{\sqrt{t}\Lambda_{1}}\dot{B}^{\frac{3}{p}+\frac{1-\theta}{1+\theta}}_{p,1})\cap
\widetilde{L}^{\frac{1+\theta}{\theta}}(T,T^{*};
(e^{\sqrt{t}\Lambda_{1}}\dot{B}^{\frac{3}{p}+\frac{\theta-1}{1+\theta}}_{p,1})}
+\|\mathcal{A}(u,\delta)\|_{\widetilde{L}^{1}(T,T^{*};e^{\sqrt{t}\Lambda_{1}}\dot{B}^{\frac{3}{p}-1}_{p,1})}
\nonumber\\
&+ \|\delta\|_{\widetilde{L}^{1}(T,T^{*};e^{\sqrt{t}\Lambda_{1}}\dot{B}^{\frac{3}{q}+2}_{q,1})}
+ \|\delta \|_{\widetilde{L}^{1+\theta}(T,T^{*};
(e^{\sqrt{t}\Lambda_{1}}\dot{B}^{\frac{3}{q}+\frac{2}{1+\theta}}_{q,1})\cap
\widetilde{L}^{\frac{1+\theta}{\theta}}(T,T^{*};
(e^{\sqrt{t}\Lambda_{1}}\dot{B}^{\frac{3}{q}+\frac{2\theta}{1+\theta}}_{q,1})}
+\|\mathcal{B}(u,\delta)\|_{\widetilde{L}^{1}(T,T^{*};
e^{\sqrt{t}\Lambda_{1}}\dot{B}^{\frac{3}{q}}_{q,1})}\nonumber\\
\leq& \|u\|_{\widetilde{L}^{1}(T,T^{*};e^{\sqrt{t}\Lambda_{1}}\dot{B}^{\frac{3}{p}+1}_{p,1})}
+\|u\|_{\widetilde{L}^{1+\theta}(T,T^{*};
(e^{\sqrt{t}\Lambda_{1}}\dot{B}^{\frac{3}{p}+\frac{1-\theta}{1+\theta}}_{p,1})\cap
\widetilde{L}^{\frac{1+\theta}{\theta}}(T,T^{*};
(e^{\sqrt{t}\Lambda_{1}}\dot{B}^{\frac{3}{p}+\frac{\theta-1}{1+\theta}}_{p,1})}
+ \|\delta\|_{\widetilde{L}^{1}(T,T^{*};e^{\sqrt{t}\Lambda_{1}}\dot{B}^{\frac{3}{q}+2}_{q,1})}
\nonumber\\
&+ \|\delta \|_{\widetilde{L}^{1+\theta}(T,T^{*};
(e^{\sqrt{t}\Lambda_{1}}\dot{B}^{\frac{3}{q}+\frac{2}{1+\theta}}_{q,1})\cap
\widetilde{L}^{\frac{1+\theta}{\theta}}(T,T^{*};
(e^{\sqrt{t}\Lambda_{1}}\dot{B}^{\frac{3}{q}+\frac{2\theta}{1+\theta}}_{q,1})}
+C\|\delta\|_{\widetilde{L}^{\infty}(T,T^{*};e^{\sqrt{t}\Lambda_{1}}\dot{B}^{\frac{3}{q}}_{q,1})}^{2}
\|\delta\|_{\widetilde{L}^{1}(T,T^{*};e^{\sqrt{t}\Lambda_{1}}\dot{B}^{\frac{3}{q}+2}_{q,1})}\nonumber\\
&+C\left(\|u\|_{\widetilde{L}^{\infty}(T,T^{*};e^{\sqrt{t}\Lambda_{1}}\dot{B}^{\frac{3}{p}-1}_{p,1})}
+ \|\delta\|_{\widetilde{L}^{\infty}(T,T^{*};e^{\sqrt{t}\Lambda_{1}}\dot{B}^{\frac{3}{q}}_{q,1})}\right)
\left(\|u\|_{\widetilde{L}^{1}(T,T^{*};e^{\sqrt{t}\Lambda_{1}}\dot{B}^{\frac{3}{p}+1}_{p,1})}
+ \|\delta\|_{\widetilde{L}^{1}(T,T^{*};e^{\sqrt{t}\Lambda_{1}}\dot{B}^{\frac{3}{q}+2}_{q,1})}\right)\nonumber\\
\leq& \frac{\zeta}{2},
\end{align}
if $T$ is sufficiently close to $T^{*}$. \eqref{eq3.31} is analogous to \eqref{eq3.25} and \eqref{eq3.26},
 which implies that the solution exists on
$[T,T^{*}]$. A contradiction with $T^{*}$ is maximal.  Thus we obtain \eqref{eq1.4}. Moreover,
 if $M_{0}\ll 1$, by using Lemma \ref{lem3.1},  we can directly choose $T=\infty$
  in \eqref{eq3.21} and \eqref{eq3.30}.
  \medskip

\textbf{\textit{Step 3:} uniqueness}

Let $(u^{1},\delta^{1})$ , $(u^{2},\delta^{2})$ be two solutions of
\eqref{eq3.17}. Denote $\widetilde{u}= u^{1}-u^{2}$,
$\widetilde{\delta}=\delta^{1}-\delta^{2}$. By using the similar
procedure of $\{(\overline{u}_{n},\overline{\delta}_{n}),
n\in\mathbb{N}\}$ is a Cauchy sequence on small interval $[0,b]$
(some small constant $b>0$), we can prove
$(\widetilde{u},\widetilde{\delta})=(0,0)$ on $[0,b]$. By repeating
the procedure, we can obtain
$(\widetilde{u},\widetilde{\delta})=(0,0)$ on $[0,T^{*})$.
 $\hfill\Box$

\section{Proof of Corollary \ref{cor1.4}}

In this section, we shall give the proof of Corollary \ref{cor1.4}. Let us first introduction the following
lemma on the operator $\Lambda^{m}e^{-\sqrt{t}\Lambda_{1}}$.

\begin{lemma}\label{lem4.1}
For all $m\geq 0$, the operator $\Lambda^{m}e^{-\sqrt{t}\Lambda_{1}}$ is a convolution operator with a kernel $k_{m}(t,\cdot)\in L^{1}(\mathbb{R}^{3})$ for all $t>0$. Moreover, we have
\begin{align*}
\|k_{m}(t,\cdot)\|_{L^{1}}\leq C_{m} t^{-\frac{m}{2}}.
\end{align*}
\end{lemma}

\begin{proof}
The proof of the Lemma is essentially due to Lemari\'{e}-Rieusset \cite{PG}, for
readers convenience, we give it as follows. By using the scaling property, it sufficient to prove
 \begin{align*}
\|k_{m}(1,\cdot)\|_{L^{1}}\leq C_{m}.
\end{align*}
 We have $\widehat{k}_{m}(1,\xi)=|\xi|^{m}e^{-|\xi|_{1}}\in L^{1}(\mathbb{R}^{3})$,
  thus $k_{m}(1,\cdot)$ is a continuous bounded function. Moreover, if $m>0$,
  we introduce a function  $\omega\in \mathcal{D}(\mathbb{R}^{3})$ so
  that $0 \notin\operatorname{Supp} \omega$ and $\sum_{j\in \mathbb{Z}} \omega(2^{j}\xi)=1$.
  Then $|\xi|^{m}\omega(\xi)\in \mathcal{D}(\mathbb{R}^{3})$,
  and if we write $|\xi|^{m}\omega(\xi)=\widehat{\Omega}_{m}(\xi)$ and
   $\theta=1-\sum_{j\geq 0}\omega(2^{j}\xi)$, we have
   $\widehat{k}_{m}(1,\xi)=\sum_{j\geq 0} 2^{-jm}\widehat{\Omega}_{m}(2^{j}\xi)e^{-|\xi|_{1}}
   +\theta(\xi)|\xi|^{m}e^{-|\xi|_{1}}$, hence
 \begin{align*}
 \|k_{m}(1,\cdot)\|_{L^{1}}\leq \sum_{j\geq 0} 2^{-jm}\|\Omega_{m}\|_{L^{1}} \|\mathcal{F}^{-1}(e^{-|\xi|_{1}})\|_{L^{1}}+\|\mathcal{F}^{-1}(\theta(\xi)|\xi|^{m}e^{-|\xi|_{1}})\|_{L^{1}}<\infty.
 \end{align*}
 Thus we complete the proof of Lemma \ref{lem4.1}.
\end{proof}
\medskip

By using Lemma \ref{lem4.1} and the result of Theorem \ref{thm1.1}, we now turn to proof of Corollary \ref{cor1.4}.
\medskip

\textbf{Proof of Corollary \ref{cor1.4}.}
 Since Theorem \ref{thm1.1} tell us that if the initial data
$(u_{0}, d_{0}-\overline{d}_{0})$ belongs to $\dot{B}^{\frac{3}{p}-1}_{p,1}(\mathbb{R}^{3})\times
\dot{B}^{\frac{3}{q}}_{q,1}(\mathbb{R}^{3})$ with $1<p,q<\infty$ satisfying \eqref{eq1.3}, then the solution
$(u,d-\overline{d}_{0})$ is in the Gevrey class $\widetilde{L}^{\infty}_{T^{*}}(
e^{\sqrt{t}\Lambda_{1}}\dot{B}^{\frac{3}{p}-1}_{p,1}(\mathbb{R}^{3}))\cap \widetilde{L}^{1}_{T^{*}}(
e^{\sqrt{t}\Lambda_{1}}\dot{B}^{\frac{3}{p}+1}_{p,1}(\mathbb{R}^{3}))
\times \widetilde{L}^{\infty}_{T^{*}}(
e^{\sqrt{t}\Lambda_{1}}\dot{B}^{\frac{3}{q}}_{q,1}(\mathbb{R}^{3}))\cap
\widetilde{L}^{1}_{T^{*}}(e^{\sqrt{t}\Lambda_{1}}\dot{B}^{\frac{3}{q}+2}_{q,1}(\mathbb{R}^{3}))$.
Consequently, for all $m\geq 0$, applying  Lemma \ref{lem4.1}, we have
\begin{align*}
&\|\Lambda_{1}^{m}u\|_{\widetilde{L}^{\infty}_{T^{*}}(\dot{B}^{\frac{3}{p}-1}_{p,1})
\cap \widetilde{L}^{1}_{T^{*}}(\dot{B}^{\frac{3}{p}+1}_{p,1})}
+\|\Lambda_{1}^{m}(d-\overline{d}_{0})\|_{\widetilde{L}^{\infty}_{T^{*}}(\dot{B}^{\frac{3}{q}}_{q,1})
\cap \widetilde{L}^{1}_{T^{*}}(\dot{B}^{\frac{3}{q}+2}_{q,1})}\nonumber\\
\lesssim& \|\Lambda_{1}^{m}e^{-\sqrt{t}\Lambda_{1}}
e^{\sqrt{t}\Lambda_{1}}u\|_{\widetilde{L}^{\infty}_{T^{*}}(\dot{B}^{\frac{3}{p}-1}_{p,1})
\cap \widetilde{L}^{1}_{T^{*}}(\dot{B}^{\frac{3}{p}+1}_{p,1})}
+\|\Lambda_{1}^{m}e^{-\sqrt{t}\Lambda_{1}} e^{\sqrt{t}\Lambda_{1}}
(d-\overline{d}_{0})\|_{\widetilde{L}^{\infty}_{T^{*}}(\dot{B}^{\frac{3}{q}}_{q,1})
\cap \widetilde{L}^{1}_{T^{*}}(\dot{B}^{\frac{3}{q}+2}_{q,1})}\nonumber\\
\lesssim& t^{-\frac{m}{2}}\left(\| e^{\sqrt{t}\Lambda_{1}}u\|_{\widetilde{L}^{\infty}_{T^{*}}(
\dot{B}^{\frac{3}{p}-1}_{p,1})\cap
 \widetilde{L}^{1}_{T^{*}}(\dot{B}^{\frac{3}{p}+1}_{p,1})}
+\|  e^{\sqrt{t}\Lambda_{1}} (d-\overline{d}_{0})\|_{\widetilde{L}^{\infty}_{T^{*}}(
\dot{B}^{\frac{3}{q}}_{q,1})\cap
\widetilde{L}^{1}_{T^{*}}(\dot{B}^{\frac{3}{q}+2}_{q,1})}\right)\nonumber\\
\lesssim& t^{-\frac{m}{2}}\left(\|u\|_{\widetilde{L}^{\infty}_{T^{*}}(
e^{\sqrt{t}\Lambda_{1}}\dot{B}^{\frac{3}{p}-1}_{p,1})\cap \widetilde{L}^{1}_{T^{*}}(
e^{\sqrt{t}\Lambda_{1}}\dot{B}^{\frac{3}{p}+1}_{p,1})}+\|(d-\overline{d}_{0})\|_{
\widetilde{L}^{\infty}_{T^{*}}( e^{\sqrt{t}\Lambda_{1}}\dot{B}^{\frac{3}{q}}_{q,1})
\cap \widetilde{L}^{1}_{T^{*}}(e^{\sqrt{t}\Lambda_{1}}\dot{B}^{\frac{3}{q}+2}_{q,1})}\right)\nonumber\\
\lesssim & t^{-\frac{m}{2}}\big(\|u_{0}\|_{\dot{B}^{\frac{3}{p}-1}_{p,1}}
+\|d_{0}-\overline{d}_{0}\|_{\dot{B}^{\frac{3}{q}}_{q,1}}\big).
\end{align*}
This completes the proof of Corollary \ref{cor1.4}.
 $\hfill\Box$
 \medskip
\\
\textbf{Acknowledgments}

The author would glad to acknowledge his sincere thanks to
Professor Song Jiang for  many valuable comments and suggestions.




\begin{thebibliography}{99}

\bibitem{BBT}H. Bae, A. Biswas, E. Tadmor, Analyticity and decay estimates of the Navier--Stokes equations in critical Besov spaces, Arch. Rational Mech. Anal., 205 (2012), 963--991.

\bibitem{BCD} H. Bahouri, J. Chemin and R. Danchin, \textit{Fourier
 Analysis and Nonlinear Partial Differential Equations}, Springer
 Heidelberg Dordrecht London New York, 2011.

\bibitem{AB} A. Biswas, Gevrey regularity for the supercritical quasi-geostrophic equation, J. Differential Equations, 257 (2014), 1753--1772.


\bibitem{CKN} L. Caffarelli, R. Kohn and L. Nireberg, Partial regularity of suitable weak solutions of Navier--Stokes
equations, Commun. Partial Differential Equations, 35 (1982), 771--831.

\bibitem{MC} M. Cannone, A generalization of a theorem by Kato on Navier--Stokes equations, Rev. Mat.
Iberoamericana, 13 (1997), 515--541.



\bibitem{DW} Y. Du and K. Wang, Space-time regularity of the Koch$\&$Tataru solutions
to the liquid crystal equations, SIAM J. Math. Anal., 45(6) (2013), 3838--3853.

\bibitem{ER} J.  Ericksen, Hydrostatic theory of liquid crystal, Arch. Rational Mech. Anal.,
 9 (1962), 371--378.


\bibitem{FT}C. Foias and R. Temam, Gevrey class regularity for the solutions of the Navier--Stokes
equations, J. Funct. Anal., 87 (1989), 359--369.

\bibitem{GPS} P. Germain, N. Pavlovic and G. Staffilani, Regularity of solutions to the Navier--Stokes
equations evolving from small data in $BMO^{-1}$, Int. Math, Res. Not., 21 (2007), 35 pages.

\bibitem{GS} Y. Giga and O. Sawada, On regularizing-decay rate estimates for solutions to the
Navier--Stokes initial value problem, Nonlinear Anal. Appl., 1(2) (2003), 549--562.

\bibitem{HL13}Y. Hao and X. Liu, The existence and blow-up criterion of liquid crystals system
in critical Besov space, Commun. Pure Appl. Anal., 13(1) (2014), 225--236.




\bibitem{JHW} J. Hineman and C. Wang, Well--posedness of nematic
liquid crystal flow in $L^{3}_{loc}(\mathbb{R}^{3})$,
Arch. Ration. Mech. Anal.,  210 (2013), 177--218.

\bibitem{HMC} M. Hong, Global existence of solutions of the
simplified Ericksen--Leslie system in dimension two, Cal. Var.
Partial Differential Equations, 40 (2011), 15--36.


\bibitem{HW} X. Hu and D. Wang, Global solution to the
three-dimensional incompressible flow of liquid crystals, Commun.
Math. Phys., 296 (2010), 861--880.

\bibitem{CHBW}C. Huang and B. Wang, Analyticity for the (generalized) Navier-Stokes equations with rough initail data, arXiv.13102.2141v2 [math.AP] 31 Oct. 2013.

\bibitem{HW1} T. Huang and C. Wang, Blow up criterion for nematic
liquid crystal flows, Commun. Partial Differential Equations, 37
(2012), 875--884.

\bibitem{TK}T. Kato, Strong $L^p$ solutions of the Navier--Stokes equations in Rm with applications to weak
solutions, Math. Z., 187 (1984), 471--480.


\bibitem{KTataru} H. Koch and  D. Tataru, Well-posedness for the Navier--Stokes equations, Adv. Math., 157 (2001),
22--35.



\bibitem{LE} F. Leslie, Theory of flow phenomenum in liquid crystals. In: The theory of
liquid crystals, London-New York: Academic Press, 4 (1979), 1--81.

\bibitem{PG}  P. G. Lemari\'{e}-Rieusset, \textit{Recent Developments in the
 Navier--Stokes Problem}, Chapman and Hall/CRC, 2002.

\bibitem{XLW} X. Li and D. Wang, Global solution to the
incompressible flow of liquid crystal, J. Differential Equations,
252 (2012), 745--767.

\bibitem{L} F. Lin, Nonlinear theory of defects in nematic liquid crystals; phase transition and flow phenomena,
 Commun. Pure Appl. Math., 42 (1989), 789--814.

\bibitem{LLW} F. Lin, J. Lin and C. Wang, Liquid crystal flow in two dimensions, Arch.
 Rational Mech. Anal., 197 (2010), 297--336

\bibitem{LL1} F. Lin and C. Liu, Nonparabolic dissipative systems modeling the flow of liquid crystals,
Commun. Pure Appl. Math.,  48  (1995), 501--537.

\bibitem{LL2} F. Lin and C. Liu, Partial regularities of the nonlinear dissipative systems modeling
the flow of liquid crystals, Discrete Contin. Dyn. Syst., A 2 (1996),
1--23.

\bibitem{LW} F. Lin and C. Wang, On the uniqueness of heat flow of harmonic maps and hydrodynamic
flow of nematic liquid crystals, Chinese Annal. Math. Ser., B 31
(2010), 921--938.


\bibitem{LW14} F. Lin and C. Wang, Global existence of weak solutions
of the nematic liquid crystal flow in dimension three, Commun. Pure Appl. Math., (2015),  Doi:10.1002/cpa/21583.



\bibitem{QL1} Q. Liu, Space-time derivative estimates of the Koch-Tataru solutios to the nematic liquid crystal system in Besov spaces,
J. Differential Equations, 258 (2015), 4368--4397.

%

\bibitem{LZZ} Q. Liu, T. Zhang and J. Zhao, Well-poseness for the 3D incompressible nematic nematic liquid crystal system in the critical  $L^{p}$ framework, Discrete Contin. Dyn. Syst., 36 (2016), 371--402.



\bibitem{MS} H. Miura and O. Sawada, On the regularizing rate estimates of Koch-Tataru's solution to
the Navier--Stokes equations, Asymptot. Anal., 49 (2006), 1--15.


\bibitem{IS} I. Stewart, \textit{The static and dynamic continuum theory of liquid crystals}, Taylor $\&$ Francis, London, New York, 2004.







\bibitem{W} C. Wang, Well-posedness for the heat flow of harmonic
maps and the liquid crystal flow with rough initial data, Arch.
Rational Mech. Anal., 200 (2011), 1--19.

\bibitem{WD} H. Wen and S. Ding, Solutions of incompressible
hydrodynamic flow of liquid crystals, Nonlinear Anal. Real World
Appl., 12 (2011), 1510--1531.

\bibitem{XZ}X. Xu and  Z. Zhang, Global regularity and uniqueness of weak solution for the 2-D liquid crystal flows,
J. Differential Equations, 252 (2012), 1169--1181.




\end{thebibliography}
\end{document}